\documentclass[11pt]{amsart2000}
\input diagrams.tex 
\newtheorem{theorem}{Theorem}[section]
\newtheorem{proposition}[theorem]{Proposition}
\newtheorem{lemma}[theorem]{Lemma}

\newtheorem{definition}[theorem]{Definition}
\newtheorem{corollary}[theorem]{Corollary}
\newcommand{\rip}[3]{\left\langle{#2},{#3}\right\rangle_{\!{#1}}}
\newcommand{\lip}[3]{{\vphantom\langle}_{#1}\!\!\left\langle{#2},
{#3}\right\rangle}
\def\radjoint{\to_{r}}
\def\ladjoint{\to_{l}}
\def\Radjoint{\to_{R}}
\def\Ladjoint{\to_{L}}
\renewcommand{\phi}{\varphi}
\newcommand{\mathify}[1]{\ifmmode{#1}\else\mbox{$#1$}\fi}
\newcommand{\set}[1]{\mathify{\{ #1 \}}}
\newcommand{\tuple}[1]{\mathify{\langle #1 \rangle}}
\newcommand{\fromto}[3]{{#1}\colon{#2}\rightarrow{#3}}
\def\hs#1#2{\left\langle #1,#2\right\rangle}
\numberwithin{equation}{section}

\begin{document}
\setlength{\unitlength}{0.01in}
\linethickness{0.01in}
\begin{center}
\begin{picture}(474,66)(0,0) 
\multiput(0,66)(1,0){40}{\line(0,-1){24}}
\multiput(43,65)(1,-1){24}{\line(0,-1){40}}
\multiput(1,39)(1,-1){40}{\line(1,0){24}}
\multiput(70,2)(1,1){24}{\line(0,1){40}}
\multiput(72,0)(1,1){24}{\line(1,0){40}}
\multiput(97,66)(1,0){40}{\line(0,-1){40}} 
\put(143,66){\makebox(0,0)[tl]{\footnotesize Proceedings of the Ninth Prague Topological Symposium}}
\put(143,50){\makebox(0,0)[tl]{\footnotesize Contributed papers from the symposium held in}}
\put(143,34){\makebox(0,0)[tl]{\footnotesize Prague, Czech Republic, August 19--25, 2001}}
\end{picture}
\end{center}
\vspace{0.25in}
\setcounter{page}{223}
\title{Morita equivalence in the context of Hilbert modules}
\author{Jan Paseka}
\thanks{Financial Support by the NATO Research Fellowship
Program and by the Ministry of Education under the project MSM 143100009
are gratefully acknowledged.}
\address{Department of Mathematics, Masaryk University\\
Jan\'a\v c{}kovo n\'a{}m. 2a\\
66295 Brno\\ 
Czech Republic}
\email{paseka@math.muni.cz}
\thanks{Jan Paseka,
{\em Morita equivalence in the context of Hilbert modules},
Proceedings of the Ninth Prague Topological Symposium, (Prague, 2001),
pp.~223--251, Topology Atlas, Toronto, 2002}
\begin{abstract} 
The Morita equivalence of m-regular involutive quantales in the context of
the theory of Hilbert $A$-modules is presented. 
The corresponding fundamental representation theorems are shown. 
We also prove that two commutative m-regular involutive quantales are
Morita equivalent if and only if they are isomorphic.
\end{abstract}
\maketitle

In the paper \cite{borceuxvi} F. Borceux and E.M. Vitale made a first step
in extending the theory of Morita equivalence to quantales. They
considered unital quantales and the category of all left modules over
these unital quantales which are unital in a natural sense. They proved
that two such module categories over unital quantales $A$ and $B$, say,
are equivalent if and only if there exist a unital $A-B$ bimodule $M$ and
a unital $B-A$ bimodule $N$ such that $M\otimes_B N\simeq A$ and
$N\otimes_A M\simeq B$.

The aim of this paper is to extend this theory in the following way:
to cover also the case of m-regular (generally non-unital) involutive
quantales and Hilbert modules over them. 
Our motivation to work in this setting comes from the theory of operator
algebras, where there is a theory of Morita equivalence for
C$^{*}$-algebras for the non-unital case (see \cite{blecher}, \cite{lance}
and \cite{RW}).
Our presentation is a combination of those in \cite{anderson}, 
\cite{blecher} and \cite{borceuxvi}.

This paper is closely related to the papers \cite{pasekahilb} and
\cite{pasekainterior} where the interested reader can find unexplained
terms and notation concerning the subject.
For facts concerning quantales and quantale modules in general we refer to
\cite{rosenthal1}. 
The algebraic background may be found in any account of Morita theory for
rings, such as \cite{anderson} or \cite{bass}.

The paper is organized as follows. 
First, we recall the notion of a right Hilbert $A$-module and related
notions. 
In Section 1 the necessary basic properties of right Hilbert modules are
established. 
Moreover, a categorical characterization of surjective module maps in the
category of m-regular right Hilbert modules is given.
In section 2 the key result is the Eilenberg-Watts theorem which states
that colimit preserving $*$-functors between categories of Hilbert modules
correspond to right Hilbert bimodules. 
Our second result is the fundamental Morita theorem for Hilbert modules.
As a consequence we get that two m-regular involutive quantales are
strongly Morita equivalent if and only if they are Morita
equivalent. 
Moreover, two commutative m-regular involutive quantales are Morita
equivalent if and only if they are isomorphic.

The present work was largely developed during the author's visit to the
Mathematics Department at the Universit\'e Catholique de Louvain, whose
members he would like to thank for their warm hospitality.

\section*{Preliminaries}

Let us begin by establishing the common symbols and notations in this
paper.

In what follows, a complete lattice will be called {\em sup-lattice}.
{\em Sup-lattice homomorphisms} are maps between sup-lattices preserving
arbitrary joins. 
We shall denote, for $S, T$ sup-lattices, $SUP(S, T)$ the sup-lattice of
all sup-lattice homomorphisms from $S$ to $T$, with the supremum given by
the pointwise ordering of mappings. 
If $S=T$ we put $\mathcal{Q}(S)=SUP(S, S)$.
Recall that a {\it quantale\/} is a sup-lattice $A$ with an associative
binary multiplication satisfying
$$
x\cdot\bigvee\limits_{i\in I} x_i = 
\bigvee\limits_{i\in I}x\cdot x_i\ \ 
\hbox{and}\ \ 
(\bigvee\limits_{i\in I}x_i)\cdot x = 
\bigvee\limits_{i\in I}x_i\cdot x
$$
for all $x,x_i\in A$, $i\in I$ ($I$ is a set).
$1$ denotes the greatest element of $A$, $0$ is the smallest element of $A$. 
A quantale $A$ is said to be {\it unital} if there is an element $e\in A$
such that $e\cdot a= a = a\cdot e$ for all $a\in A$. 
An {\em opposite quantale} $A^{d}$ to $A$ is a sup-lattice with the same
join operation as $A$ and with the multiplication $a\circ b=b\cdot a$.
A {\it subquantale\/} $A'$ of a quantale $A$ is a subset of $A$ closed
under $\bigvee$ and $\cdot$.
Since the operators $a\cdot -$ and $-\cdot b: A\to A$, $a, b\in A$
preserve arbitrary joins, they have right adjoints. 
Explicitly, they are given by
$$
a\to_{r}c = \bigvee\{s\in A| a\cdot s\leq c\}
\quad
\hbox{and}
\quad 
b\to_{l}d = \bigvee\{t\in A| t\cdot b\leq d\}
$$
\noindent respectively.

An {\em involution} on a sup-lattice $S$ is a unary operation such that
$$
\begin{array}{rcl}
a^{**}&
=
&a,\\
(\bigvee a_{i})^{*}&
=&
\bigvee a_{i}^{*}
\end{array}
$$
for all $a, a_{i}\in S$.
An {\em involution} on a quantale $A$ is an involution on the sup-lattice 
$A$ such that
$$
\begin{array}{rcl}
(a\cdot b)^{*}&
=&
b^{*}\cdot a^{*}\\
\end{array}
$$
for all $a, b\in A$.
A sup-lattice (quantale) with the involution is said to be {\em
involutive}.

By a {\em morphism of} ({\em involutive}) {\em quantales}
will be meant a $\bigvee$- ($^{*}$-) and {$\cdot$-preserving} mapping 
$f:A\to A'$. 
If a morphism preserves the unital element we say that it is {\em unital}.

Let $A$ be a quantale.
A {\em right module over} $A$ (shortly a right $A$-module) is a
sup-lattice $M$, together with a {\em module action} 
$${-}{\diamond}{-}:M\times A\to M$$ 
satisfying
\begin{itemize}
\item[(M1)]
$m\diamond(a\cdot b) = (m\diamond a)\diamond b$
\item[(M2)]
$(\bigvee X)\diamond a = \bigvee \{x\diamond a: x\in X\}$
\item[(M3)]
$m\diamond \bigvee S = \bigvee \{m\diamond s: s\in S\}$
\end{itemize}
for all $a, b\in A$, $m\in M$, $S\subseteq A$, $X\subseteq M$.
So we have two maps 
$$
- \Ladjoint -:M\times A\to M,
\quad
- \Radjoint -:M\times M\to A
$$ 
such that, for all $a\in A$, $m,n \in M$, 
$$
m\diamond a\leq n
\quad
\mbox{iff}
\quad
a\leq m\Radjoint n
\quad
\mbox{iff} 
\quad 
m\leq a\Ladjoint n.
$$
$M$ is called a {\em unital} $A$-module if $A$ is a unital
quantale with the unit $e$ and $m\diamond e=m$ for all $m\in M$.

Let $M$ and $N$ be modules over $A$ and let $f:M\to N$ be a sup-lattice
homomorphism. $f$ is a {\em module homomorphism} if 
$f(m\diamond a)=f(m)\diamond a$ for all $a\in A, m\in M$.
We shall denote by $MOD_A$ the category of right $A-$modules and module
homomorphisms.

For a module $X$ in $MOD_A$ the submodule $\hbox{\rm ess}(X)=X\diamond A$
generated by the elements $x\diamond a$ is called the {\em essential part}
of $X$.
If $\hbox{\rm ess}(X)=X$ we say that $X$ is {\em essential}.
The full subcategory of essential $A$-modules is denoted 
$\hbox{\rm ess}-MOD_A$. 
We shall say that $A$ is {\em right separating} for the $A$-module $M$
and that $M$ is ({\em right}) {\em separated} by $A$ if 
$m\diamond (-)=n\diamond (-)$ implies $m=n$.
We say that $M$ is {\it m-regular} if it is both separated by $A$ and
essential. 
An involutive quantale $A$ is called {\it m-regular} if it is {m-regular}
as an $A$-module.
Then, evidently $1\cdot 1=1$ in $A$ and $a\cdot (-)=b\cdot (-)$ implies
$a=b$.

Note that we may dually define the notion of a (unital) left $A$-module
with a left multiplication $\bullet$. 
We then have two maps 
$$
- \ladjoint -:M\times M\to A,
\quad 
- \radjoint -:A\times M\to M
$$
such that, for all $a\in A$, $m,n \in M$, 
$$
a\bullet m\leq n
\quad
\mbox{iff} 
\quad
a\leq m\ladjoint n
\quad
\mbox{iff} 
\quad
m\leq a\radjoint n.$$

The theory of Hilbert $A$-modules (we refer the reader to 
\cite{pasekahilb} for details and examples) is a generalization of the
theory of complete semilattices with a duality and it is the natural
framework for the study of modules over an involutive quantale $A$ endowed
with $A$-valued inner products.

Let $A$ be an involutive quantale,
$M$ a right (left) $A$-module. We say that $M$ is a 
{{\em right 
\hbox{\rm{}(}left\hbox{\rm{})}} Hilbert $A$-module}
{\hbox{\rm{}(}\em right \hbox{\rm{}(}left\hbox{\rm{})} strict Hilbert
$A$-module),
{{\em right \hbox{\rm{}(}left\hbox{\rm{})}} pre-Hilbert $A$-module}}
if $M$ is equipped with a map 
$$\langle -, - \rangle:M\times M\to A$$
called the {\em inner product}, such that for all $a\in A$, $m, n\in M$
and $m_i\in M$, where $i\in I$, the conditions
{\rm{}(\ref{whspi})--(\ref{whspiv})
((\ref{whspi})--(\ref{whspvi}), 
(\ref{whspi})--(\ref{whspv}))}
are satisfied.

\begin{eqnarray}
&
\quad&
\langle m, n\rangle\cdot a=\langle m, n\diamond a\rangle\quad
\hbox{\rm{}(}a\cdot \langle m, n\rangle=\langle a\bullet m, n\rangle\hbox{\rm{})};\label{whspi}\\
&\quad&\bigvee_{i\in I}
\langle m_i, n\rangle=
\langle\bigvee_{i\in I}m_i, n\rangle;\label{whspii}\\
&\quad&\bigvee_{i\in I}
\langle m ,m_i\rangle=\langle m, 
\bigvee_{i\in I}m_i\rangle;\label{whspiii}\\
&\quad&
\langle m, n\rangle^{*}=\langle n, m\rangle;\label{whspv}\\
&\quad&
\langle - ,m\rangle=\langle -, n\rangle\
(\langle m , -\rangle=\langle n, -\rangle) \
\hbox{implies}\ m=n;\label{whspiv}\\
&\quad&
\langle m ,m\rangle=0 \
\hbox{implies}\ m=0.\label{whspvi}
\end{eqnarray}

If $A$ is an involutive quantale, let $HMOD_A$ be the category of right
Hilbert $A$-modules with morphisms the usual $A$-module maps. 
The full subcategory of m-regular right Hilbert $A$-modules is denoted
$\hbox{\rm mreg}-HMOD_A$.

Let $A$ be an involutive quantale, $f:M\to N$ a map between right
(left) pre-Hilbert $A$-modules.
We say that a map $g:N\to M$ is a {\em $^{*}$-adjoint to $f$} and
$f$ is {\em adjointable} if
$$
\langle f(m), n\rangle = \langle m, g(n)\rangle
$$
for all $m\in M$, $n\in N$. 
Evidently, any adjointable map is a module homomorphism. 
If $f$ is adjointable we put
$$
f^{*} = 
\bigvee\{g:N\to M; g\ \hbox{is a $^{*}$-adjoint to $f$}\}.
$$

Note that $f\leq f^{**}$ and $f^{*}=f^{***}$. 
If $M$ and $N$ are Hilbert $A$-modules the $^{*}$-adjoint to $f$ is
uniquely determined by property (\ref{whspiv}) i.e. $f=f^{**}$.
The set of all adjointable maps from $M$ to $N$ is denoted by 
$\mathcal{A}_{A}(M,N)$. 
We shall denote by $\hbox{\rm mreg}-Hilb_A$ the category of m-regular
right Hilbert $A$-modules and adjointable mappings. 
We say that an adjointable map $f:M\to N$ is an {\em isometry} if, for all
$m_1, m_2\in M$, $\hs{m_1}{m_2}=\hs{f(m_1)}{f(m_2)}$. 
This is equivalent to $f^{*}\circ f=\hbox{\rm id}_M$. 
Similarly, an adjointable map $f:M\to N$ is {\em unitary} if 
$f^{*}\circ f=\hbox{\rm id}_M$ and $f\circ f^{*}=\hbox{\rm id}_N$. 
Note that any surjective isometry is necessarily unitary.

Recall that from \cite{pasekahilb} we know that, for any right Hilbert
$A$-module $M$ and for all $m\in M$, the map $m^{\sim}:A\to M$ defined by
$a\mapsto m\diamond a$ has a $^{*}$-adjoint $m^{\star}:M\to A$ defined by
$n\mapsto \langle m, n\rangle$.

Similarly, let $A$ and $B$ be involutive quantales, and let $M$ and $X$ be
right Hilbert $B$-modules.
We say that $M$ is a {\em right Hilbert $A-B$ bimodule\/} if it is
a left $A$-module satisfying
\begin{equation}\label{rHb-eq}
a\bullet (x\diamond b) = (a\bullet x)\diamond b\
\mbox{and}\
\rip{B}{a\bullet x}{y} = \rip{B}{x}{a^*\bullet y}
\end{equation}
for all $a\in A$, $x,y\in M$, and $b\in B$.
We say that $F$ is an m-regular right Hilbert $A-B$ bimodule if it is
both an m-regular left $A$-module and an m-regular right $B$-module.
An {\em isomorphism} of right Hilbert bimodules is a bijective
$\bigvee$-preserving map $\Phi\colon M \to N$ such that
\begin{itemize}
\item[(i)] $\Phi(a\bullet x) = a\bullet\Phi(x)$,
\item[(ii)] $\Phi(x\diamond b) = \Phi(x)\diamond b$, and
\item[(iii)] $\rip{B}{\Phi(x)}{\Phi(y)} = \rip{B}{x}{y}$
\end{itemize}
i.e.\ $\Phi$ is an isometric surjective bimodule homomorphism. 
In particular, $\Phi$ is a unitary $B$-module map.

We say that a (full m-regular) right Hilbert $A-B$ bimodule $X$ is an
({\em imprimitivity}) $A-B$ {\em bimodule} if
$X$ is also a (full m-regular) left Hilbert $A$-module in such a way that
$$
\lip{A}{x}{y}\bullet z = x\diamond\hs{y}{z}_{B}\ 
\mbox{and}\
\lip{A}{x\diamond b}{y} = \lip{A}{x}{y\diamond b^*}.
$$

Suppose that $M$ is a right Hilbert $A$-module, $N$ is a right Hilbert
$A-B$ bimodule. 
Then the sup-semilattice tensor product $M\otimes N$ is a pre-Hilbert
$B$-module under the pre-inner product given on simple tensors by
$$
\hs{x_1\otimes y_1}{x_2\otimes y_2}_{M\otimes N} =
\hs{y_1}{\langle x_1, x_2\rangle\bullet y_2}
$$
for all $x_1, x_2\in M$ and $y_1, y_2\in N$. 
We then denote by $M\dot{\otimes}_{A}N$ the corresponding Hilbert
$A$-module $(M\otimes N)_{R_H}$ and we say that $M\dot{\otimes}_{A}N$ is a
{\em interior tensor product} of $M$ and $N$ over $A$.
It is shown in \cite{pasekainterior} that our interior tensor product has
similar properties as its Hilbert C$^{*}$-module counterpart. 
Note only that a right Hilbert $A$-module is m-regular if and only if 
$M\dot{\otimes}_{A}A\simeq M$ via the standard isomorphism
$m\dot{\otimes}_{A}a\mapsto m\diamond a$.

We shall sometimes use ${\bf C}_A$ for some of the following
categories: 
$\hbox{\rm mreg}-HMOD_A$,
$\hbox{\rm mreg}-Hilb_A$.
Similarly, we have the categories of left modules
${}_A{}MOD$, ${}_A{}HMOD$,
$\hbox{\rm mreg}-{}_A{}HMOD$,
$\hbox{\rm mreg}-{}{}_A{}Hilb$).

In this paper we are concerned with functors between categories of 
modules. 
Such functors, e.g.\
$F: \hbox{\rm mreg}-{Hilb}_A \rightarrow \hbox{\rm mreg}-{Hilb}_B$ 
are assumed to be join-preserving (SUP-functors) on
$\bigvee$-semilattices of morphisms. 
Thus the map $T \mapsto F(T)$ from $\hbox{\rm mreg}-{Hilb}_A(X,W)$
to $\hbox{\rm mreg}-{Hilb}_B(F(X),F(W))$ is join-preserving, for all pairs
of objects $X,W \in $ ${\bf C}_A$.

In what follows let us assume that $A$ and $B$ are m-regular involutive
quantales.
We say that a SUP-functor 
$$F : \hbox{\rm mreg}-{Hilb}_A \rightarrow \hbox{\rm mreg}-{Hilb}_B$$
is a {\em $*$-functor} if $F(T^*) = F(T)^*$ for all adjointable $A$-module
maps $T$. 
In particular, from \cite{pasekainterior} we know that the function
$$(-)\dot{\otimes}_A N : 
\hbox{\rm mreg}-{Hilb}_A\to \hbox{\rm mreg}-{Hilb}_B$$
which assigns to each right Hilbert $A$-module $M$ the inner tensor
product $M\dot{\otimes}_A N$ and to each adjointable map $f$ between right
Hilbert $A$-modules the adjointable map $f\dot{\otimes}_A \hbox{\rm id}_F$
between right Hilbert $B$-modules is a *-functor preserving biproducts.

We say two $*$-functors
$$F_1, F_2 : \hbox{\rm mreg}-{Hilb}_A \rightarrow \hbox{\rm 
mreg}-{Hilb}_B$$
are (naturally) unitarilly isomorphic, if they are naturally isomorphic
via a natural transformation $\tau$ in the sense of category theory 
\cite{anderson}, with the natural transformations being unitaries i.e.\
$\tau^{*}_{F_1(M)}\circ \tau_{F_1(M)}=\hbox{\rm id$_{F(M_1)}$}$
and
$\tau_{F(M_1)}\circ \tau^{*}_{F(M_1)}=\hbox{\rm id$_{F(M_2)}$}$.
In this case we write $F_1 \cong F_2$ {\em unitarilly}. 
Similarly, we say that a $*$-functor 
$F: \hbox{\rm mreg}-{Hilb}_A \rightarrow \hbox{\rm mreg}-{Hilb}_B$
is a {\em unitary equivalence functor} if there is a $*$-functor 
$G: \hbox{\rm mreg}-{Hilb}_B \rightarrow \hbox{\rm mreg}-{Hilb}_A$
such that we have natural unitary isomorphisms
$\eta:GF\to \hbox{\rm Id}_{\hbox{\tiny \rm mreg}-{Hilb}_A}$ 
and
$\zeta:FG\to \hbox{\rm Id}_{\hbox{\tiny \rm mreg}-{Hilb}_B}$.

\section{Hilbert modules}

In this section we develop the basic properties of Hilbert modules not
mentioned in \cite{pasekahilb}.
Let us begin with a small observation about Hilbert modules.

\begin{lemma}\label{misep} 
For any involutive quantale $A\in HMOD_A$ and any right Hilbert $A$-module
$M$, $A$ is right separating for $M$.
\end{lemma}

\begin{proof} 
Assume that $m\diamond (-)=n\diamond (-)$ i.e.\ $m\diamond a=n\diamond a$
for all $a\in A$. 
Hence, for all $p\in M$ and all $a\in A$ we have that 
$\hs{p}{m\diamond a}=\hs{p}{n\diamond a}$ 
i.e.\
$\hs{p}{m}\cdot a=\hs{p}{n}\cdot a$ 
i.e.\
$\hs{p}{m}=\hs{p}{n}$. 
So we get that $m=n$.\
\end{proof}

The preceding observation then yields

\begin{corollary}\label{epis}
Let $A$ be an involutive quantale, $A\in HMOD_{A}$.
Then a morphism $f$ in 
${\bf C}_{A}\in \{ \hbox{\rm mreg}-HMOD_A, \hbox{\rm mreg}-Hilb_A\}$
is mono iff the map $f$ is one-to-one.
\end{corollary}

\begin{proof}
Evidently, any injective morphism in ${\bf C_A}$ is a monomorphism. 
Conversely, let $f:M\to N$ be a monomorphism in ${\bf C}_A$, 
$f(m_1)=f(m_2)$. 
Then 
$f\circ (m_1\diamond (-))=f\circ (m_2\diamond (-))$ 
i.e.\
$m_1\diamond (-)=m_2\diamond (-)$. 
Hence $m_1=m_2$ and $f$ is injective.
\end{proof}

Note that in the category $Hilb_A$ over $A\in HMOD_{A}$ an adjointable map
is mono iff $f^{*}$ is epi. 
This immediately yields that for any surjective adjointable map
its adjoint is one-to-one.

\begin{lemma}\label{prodessH} 
The category 
${\bf C}_{A}\in\{\hbox{\rm mreg}-HMOD_A, \hbox{\rm mreg}-Hilb_A\}$ 
has biproducts.
\end{lemma}

\begin{proof}
We know from \cite{pasekahilb} that $Hilb_{A}$ has biproducts (with
injections $i_j$ and projections $\pi_j$) and these are exactly cartesian
products. 
The lemma follows then from the fact that evidently any cartesian product
of m-regular modules is m-regular.
\end{proof}

As usual, for any $A\in HMOD_{A}$, $A^J$ will be viewed as a Hilbert
$A$-module equipped with the inner-product
$$\langle (a_j)_{j\in J}, (b_j)_{j\in J}\rangle = 
\bigvee_{j\in J} a_j^* b_j.$$

Let us observe that $A^J$ will often stand for the set of column matrices
over $A$ of the type $S\times 1$.
In that way, the above inner-product can be expressed, for 
$v=(a_j)_{j\in J}$ and $w=(b_j)_{j\in J}$, as $\langle v, w\rangle =
v^*w$.
Note that $v^*$ refers to the conjugate-transpose matrix.
We shall denote, for any $a\in A$ and all $j\in J$, by ${\bf a}_{j}$ an
element of $A^{J}$ defined as follows:
$$
\pi_k({\bf a}_j)
=
\left\{\begin{array}{ll}
a&\hbox{\rm if}\ k=j\\
0&\hbox{otherwise.}
\end{array}\right.
$$

For each $J$-tuple $\mu = (m_j)_{j\in J}$ in the Hilbert $A$-module $M^J$, 
we denote by $\Omega_\mu$ the operator in $\mathcal{A}_{A}(A^J, M)$ defined 
by
$$
\Omega_\mu\bigl((a_j)_{j\in J}\bigr) = \bigvee_{j\in J} m_j\diamond a_j, 
\quad
(a_j)_{j\in J} \in A^J.
$$
It is easy to see that $\Omega_\mu^*$ is given by
$$
\Omega_\mu^*(x) = \bigl(\langle m_j, x\rangle\bigr)_j, 
\quad 
x \in M.$$

If $\nu = (n_j)_{j\in J}$ is an $J$-tuple of elements of $N$, then the 
operator $T = \Omega_\nu\circ \Omega_\mu^*$ is in $\mathcal{A}_A(M, N)$.
In particular, for all $m\in M$ and $n\in N$,
$\Omega_{(n)} = n^{\sim}$, $\Omega_{(m)}^{*}=m^{\star}$ and
defining 
$$\Theta_{n, m}(x) = 
(\Omega_{(n)}\circ \Omega^{*}_{(m)})(x) = n\diamond \langle m,x\langle$$
we have
$$
T(x) = \bigvee_{j\in J} \Theta_{n_j, m_j}(x), 
\quad x \in M.
$$

Maps such as $T$ are exactly {\it compact\/} operators in the sense of
\cite{pasekahilb} and the set of all compact maps will be denoted 
$\mathcal{K}_{A}(M, N)$ or just $\mathcal{K}_A(M)$ in case $M = N$. 
Note that a composition of a compact operator with an adjointable one is
again compact i.e.\
$\mathcal{K}_{A}(N, P)\circ \mathcal{A}_{A}(M, N) \subseteq 
\mathcal{K}_{A}(M, P)$ 
and
$\mathcal{A}_{A}(N, P)\circ \mathcal{K}_{A}(M, N) \subseteq 
\mathcal{K}_{A}(M, P)$.
An expository treatment of compact operators on Hilbert $A$-modules may be
found in \cite{pasekahilb}.

\begin{lemma}\label{esscom} 
Let $A$ be an m-regular involutive quantale and let 
$M\in \hbox{\rm mreg}-HMOD_A$, $M$ is full.
Then $\mathcal{K}_{A}(M)$ is an m-regular involutive quantale and $M$ is an
m-regular full left
Hilbert $\mathcal{K}_{A}(M)$-module.
\end{lemma}

\begin{proof} 
Evidently, $\mathcal{K}_{A}(M)$ is an involutive quantale (see
\cite{pasekahilb}).
Let $f, g\in \mathcal{K}_{A}(M)$.
Assume that $f\circ \Theta_{y, x}=g\circ \Theta_{y, x}$ for all
$y, x\in X$. 
This gives us $\Theta_{f(y), x}=\Theta_{g(y), x}$ and since $M$ is full
then $f(y)\diamond a=g(y)\diamond a$ for all $a\in A$ i.e. $f=g$. 
Similarly,
$$
\begin{array}{lll}
\Theta_{y, x}&
=&
\Theta_{\bigvee_{i\in I}y_i\diamond a_i, x}\\
&
=&
\bigvee_{i\in I}\Theta_{y_i\diamond a_i, x}\\
&
=&
\bigvee_{i\in I}
\Theta_{y_i\diamond \bigvee_{j\in J}\langle v_{ij}, u_{ij}\rangle, x}\\
&
=&
\bigvee_{i\in I, j\in J}
\Theta_{y_i\diamond \langle v_{ij}, u_{ij}\rangle, x}\\
&
=&
\bigvee_{i\in I, j\in J}
\Theta_{y_i, v_{ij}} \circ \Theta_{u_{ij}, x}
\end{array}
$$
for suitable elements $y_i, u_{ij}, v_{ij}$, $i\in I, j\in J$.
This gives us that $\mathcal{K}_{A}({}M)$ is m-regular.

We shall define the module action on 
$\bullet:\mathcal{K}_{A}(M)\times M\to M$ by 
$$f\bullet m=f(m)$$ 
and the inner product
$\lip{\mathcal{K}_{A}(M)}{-}{-}:M\times M\to \mathcal{K}_{A}(M)$ 
by
$$\lip{\mathcal{K}_{A}(M)}{y}{x}=\Theta_{y, x}.$$ 
Then evidently the pre-Hilbert module conditions are satisfied,
$$
M = \mathcal{K}_{A}(M)\bullet M,\
\lip{\mathcal{K}_{A}(M)}{M}{M} = \mathcal{K}_{A}(M).
$$
Assume that
$\lip{\mathcal{K}_{A}(M)}{y_1}{-} = \lip{\mathcal{K}_{A}(M)}{y_2}{-}$ 
i.e.\ $y_1\diamond a=y_{2}\diamond a$ for all $a\in A$.
Hence $y_1=y_2$.
\end{proof}

\begin{proposition}\label{cahrc} 
Let $A$ be an m-regular involutive quantale and let 
$M\in \hbox{\rm mreg}-HMOD_A$.
For each $\mu=(\mu_j)_{j\in J}$ in $M^J$ one has that $\Omega_\mu$ is in
$\mathcal{K}_A(A^J, M)$ and hence also that $\Omega_\mu^*$ is in 
$\mathcal{K}_A(M, A^J)$. 
Moreover, $\mathcal{K}_A(A^J, M)$ $\simeq$ $M^{J}$ as sup-lattices.
\end{proposition}

\begin{proof}
It is obviously enough to consider the case $|J|=1$. 
Let $\mu=(m)$ and 
$\bigvee_{\lambda\in\Lambda} m_\lambda\diamond a_{\lambda} = m$.
Therefore we have for all $a$ in $A$
$$
\begin{array}{lll}
\Omega_\mu(a)&
=& 
m\diamond a\\ 
&
=& 
(\bigvee_{\lambda\in\Lambda} m_\lambda\diamond a_{\lambda})\diamond a\\
&
=&
\Omega_{(m_\lambda)_{\lambda\in\Lambda}}
\Omega_{({a_{\lambda}}^{*})_{\lambda\in\Lambda}}^{*}(a).
\end{array}
$$

We have that
$$
\Omega_{(m)}\circ \Omega_{(a)}^{*}(b) = 
\Omega_{(m\diamond a^{*})}(b)
$$
i.e.\ 
any generator of $\mathcal{K}_A(A^J, M)$ has the form $\Omega_{(m)}$. 
Since $\Omega_{(m)}$ is compact for all $m$
we have that
$$\mathcal{K}_A(A^J, M)\simeq \{\Omega_{(m)}: m\in M\}\simeq M.$$
\end{proof}

\begin{lemma}\label{crxcahrc} 
Let $A$ be a unital involutive quantale and let 
$M\in \hbox{\rm mreg}-HMOD_A$. 
Then $HMOD_A(A^{J},M)=\mathcal{K}_A(A^J, M)\simeq M$.
\end{lemma}

\begin{proof}
Let $f:A^{J}\to M$ be any module map, $x=(x_j)_{j\in J}$. 
Then 
$$
f(x) 
=
\bigvee_{j\in J} f({\bf e}_{j})\diamond x_j
=
\bigvee_{j\in J}
\left(
\Omega_{f({\bf e}_{j})}\circ i_j\circ \Omega_{{\bf e}_{j}}^{*}
\right)(x)
.
$$
Hence,
$$
f 
=
\bigvee_{j\in J}
\Omega_{f({\bf e}_{j})}\circ i_j\circ \Omega_{{\bf e}_{j}}^{*}
\in 
\mathcal{K}_A(A^J, M)
.
$$
\end{proof}

\begin{corollary}\label{xcahrc} 
Let $A$ be a unital involutive quantale. 
Then $\hbox{\rm id}_{A^{J}}$ is in $\mathcal{K}_A(A^J)$.
\end{corollary}

\begin{definition} 
Let $A$ be an involutive quantale. 
A Hilbert $A$-module $M$ will be said to be a {\em nuclear} module if the
identity operator $\hbox{\rm id}_M$ is in $\mathcal{K}_A(M)$.
\end{definition}

Note that any m-regular involutive quantale $A\simeq\mathcal{K}_{A}(A)$ is
nuclear if and only if it is unital.
We shall now give the complete characterization of nuclear Hilbert modules
over unital involutive quantales.

\begin{proposition} 
Let $A$ be a unital involutive quantale and let 
$M\in \hbox{\rm ess}-HMOD_A$. 
The following conditions are equivalent:
\begin{enumerate}
\item $M$ is nuclear.
\item $M$ is a retract of $A^{J}$ in $Hilb_A$.
\end{enumerate}
\end{proposition}

\begin{proof}
(1)$\implies$(2).
Assume $M$ to be nuclear.
Then $\hbox{\rm id}_M = \Omega_\nu \Omega_\mu^*$ where $\mu$ and $\nu$ are
in $M^J$.

(2)$\implies$(1). 
Let $M$ be a retract of $A^{J}$ in $Hilb_A$, $r:A^{J}\to M$ the retraction
and $i:M\to A^{J}$ the embedding such that $r\circ i=\hbox{\rm id}_M$.
Since $\hbox{\rm id}_{A^{J}}\in \mathcal{K}_{A}(A^{J})$ we have that
$\hbox{\rm id}_M 
= r\circ \hbox{\rm id}_{A^{J}}\circ i
\in \mathcal{K}_{A}(M)$.
\end{proof}

Sometimes we shall need the following lemma, a part of which is contained
in \cite{pasekahilb}.

\begin{lemma}\label{pregihi} 
Let $A$ be an involutive quantale and let $M$ be a left
(right) pre-Hilbert $A$-module. 
Then the factor module $M_{R_{H}}$ defined by the equivalence relation
$$
R_{H}=\{(m, n)\in M\times M:
\hs{m}{p}=\hs{n}{p}\ \hbox{\rm for all }\ p\in M \}
$$ 
is a left (right) Hilbert $A$-module. 
Moreover, $j_H:M\to M_{R_H}$ is an adjointable map and if $f:N\to M$ is a
module (adjointable) map then there is a unique module (adjointable) map
$\overline{f}:N\to M_{R_H}$ such that $\overline{f}=j_{H}\circ f$; here
$j_{H}(m)=\bigvee\{n\in M:(m, n)\in R_{H}\}$.
Similarly, if $g:M\to P$ is an adjointable map then there is a unique 
adjointable map $\hat{g}:M_{R_H}\to P$ such that $\hat{g}=g\circ
{j_{H}}^{*}$.
\end{lemma}

\begin{proof} 
The main part of this lemma was proved in \cite{pasekahilb} and the
adjointability of $j_H$ follows from its definition. 
So let $f:N\to M$ be a module (adjointable) map.
We shall put $\overline{f}(n)=j_H(f(n))$. 
Evidently, $\overline{f}$ is a module map.
Let us check that $\overline{f}$ is adjointable. 
Assume $n\in N, j_{H}(m)\in M_{R_H}$. 
Then
$$
\hs{\overline{f}(n)}{j_H(m)} 
= \hs{f(n)}{j_H(m)} 
= \hs{n}{f^{*}(j_H(m))}.
$$ 
Clearly, such $\overline{f}$ is uniquely determined.

Similarly, let $g:M\to P$ be an adjointable map. 
We define 
$$\hat{g}(m)=g({j_H}^{*}(m))=g(m).$$
Since $\hat{g}$ is a composition of adjointable maps it is adjointable.
The uniqueness is evident.
\end{proof}

\begin{theorem}\label{prehilblim} 
Let $A$ be an involutive quantale. 
Then the category of pre-Hilbert $A$-modules has limits of arbitrary
diagrams.
\end{theorem}

\begin{proof}
The proof follows general category theoretic principles. 
We describe the limit of the diagram 
$$\tuple{(M_i)_{i\in O},(\fromto{f_j}{M_{d(j)}}{M_{c(j)}})_{j\in J}}$$ 
as a set of particular elements of the product of all $M_i$'s, the
so-called {\em commuting tuples}.
$$
M = 
\set{(x_i)_{i\in O}\in \prod_{i\in O}M_i : 
\forall j\in J\ x_{c(j)}=f_j(x_{d(j)})}
.
$$
Evidently, $M$ is a $A$-submodule of the product, that is, $M$ is a
pre-Hilbert $A$-module because the coordinatewise supremum of commuting
tuples is commuting as all $f_j$ are module maps. 
This also proves that the projections 
$\fromto{\pi_j}{\prod_{i\in O}M_i}{M_j}$ restricted to $M$ are module
maps. 
They give us the maps needed to complement $M$ to a cone.

Given any other cone $\tuple{E,(\fromto{g_i}E{M_i})_{i\in O}}$, we define
the mediating morphism $\fromto{h}{E}{M}$ by 
$h(x)=(g_i(x))_{i\in O}$. 
Again, it is obvious that this is well-defined and a module map, and that
it is the only possible choice.
\end{proof}

We also have the dual:

\begin{theorem}\label{colim}
Let $A$ be an involutive quantale.
Then the category of pre-Hilbert $A$-modules has colimits of arbitrary
diagrams.
\end{theorem}

\begin{definition}\label{faicogen}
Let {\bf C} be a subcategory of $MOD_A$ and let $M$ be a module in {\bf C}.
\begin{enumerate}
\item 
$M$ is called {\rm faithful} if
$m\diamond a=m\diamond b$ for all $m\in M$ implies that $a=b$.
\item 
$M$ {\em generates} $X\in$~{\bf C} if there is a direct sum
$\coprod_{\gamma\in \Gamma} M$ of copies of $M$ and an epimorphism
$\phi:\coprod_{\gamma\in \Gamma} M\to X$ in {\bf C}.
$M$ is a {\em generator} for {\bf C}, if $M$ generates all modules
$X\in$~{\bf C}.
\item 
$M$ {\em cogenerates} $X\in$~{\bf C} if there is a direct product
$\prod_{\gamma\in \Gamma} M$ of copies of $M$ and a monomorphism
$\psi:X\to\prod_{\gamma\in \Gamma} M$ in {\bf C}.
$M$ is a {\em cogenerator} for {\bf C}, if $M$ cogenerates all modules
$X\in$~{\bf C}.
\end{enumerate}
\end{definition}

\begin{lemma}\label{fulfait}
Let $A$ be a right separating involutive quantale, $M$ be a full right
Hilbert $A$-module. 
Then $M$ is faithful.
\end{lemma}

\begin{proof}
Let $a, b\in A$ such that
$m\diamond a=m\diamond b$ for all $m\in M$. 
Then
$\hs{n}{m\diamond a}=\hs{n}{m\diamond b}$ for all $m, n\in M$ 
i.e.\
$\hs{n}{m}\cdot a=\hs{n}{m}\cdot b$ for all $m, n\in M$ 
i.e.\
$a=b$.
\end{proof}

\begin{lemma}\label{cogeneration} 
Let $A$ be an m-regular involutive quantale.
Then $A$ is both a generator and a cogenerator for 
${\bf C}_{A}\in \{ \hbox{\rm mreg}-HMOD_A, \hbox{\rm mreg}-Hilb_A\}$.
More exactly, for every $M$ in ${\bf C}_{A}$, $M$ is both a quotient
module and a submodule of the sum 
$\bigsqcup_{m\in M} A\cong \prod _{m\in M} A=A^{M}$.
\end{lemma}

\begin{proof} 
We shall use the compact maps
$m^{\sim}:A\to M$ and $m^{\star}:M\to A$, $m\in M$.

$$
\begin{diagram}
A^{M}&\rTo^{{\bf p}_M}&M&\rTo^{{\bf i}_M}&A^{M}\\
&\luTo_{i_m} &\uTo^{m^{\sim}}\dTo_{m^{\star}}&
\ldTo_{\pi_m} & \\
& &A& &\\
\end{diagram}
$$

They give us compact maps 
${\bf p}_{M}:A^{M}\to M$, 
${\bf p}_{M}=\Omega_{(m)_{m\in M}}$
and 
${\bf i}_M:M\to A^{M}$, 
${\bf i}_{M}=\Omega_{(m)_{m\in M}}^{*}$;
here 
$\pi_m\circ {\bf i}_M=m^{\star}$ and 
${\bf p}_{M}\circ {\iota}_m=m^{\sim}$.
Since $M$ is a Hilbert $A$-module ${\bf i}_{M}$ is injective and since $M$
is essential ${\bf p}_{M}$ is surjective. 
Namely, for any $n\in M$, there is an element $u_n\in A^M$ defined by
$$
p_{m}(u_{n})=\bigvee\{a\in A: m\diamond a\leq n\}
$$
such that
$n=\bigvee_{m\in M} m\diamond {p_m}(u_n)={\bf p}_M(u_n)$.
\end{proof}

Note that the surjective module map ${\bf p}_{M}:A^{M}\to M$ may be
defined for any essential $A$-module $M$.
Similarly as in \cite{gronbak} we have the following proposition.

\begin{proposition} 
Let $A$ be an m-regular involutive quantale. 
Then $U\in \hbox{\rm mreg}-HMOD_{A}$
($U\in \hbox{\rm mreg}-Hilb_{A}$) is faithful if and only if it
cogenerates a generator.
\end{proposition}

\begin{proof}
Suppose that $U$ is faithful.
Then the module homomorphism (adjointable map) $f:A\to U^{U}$ defined by
$f(a)(u)=u\diamond a$ is monic and $U$ cogenerates the generator $A$. 
Conversely, let $M$ be a generator for $\hbox{\rm mreg}-HMOD_{A}$ and $U$
a cogenerator of $M$. 
Let $u\diamond a=u\diamond b$ for all $u\in U$. 
Then $m\diamond a=m\diamond b$ since $M$ is embeddable into the product of
copies of $U$. 
Since we have an epimorphism $g:M^{J}\to A$ such that 
$a^{\sim}\circ g=b^{\sim}\circ g$ we obtain that $c\cdot a=c\cdot b$ for
all $c\in A$ i.e.\ $a=b$.
\end{proof}

\begin{proposition}\label{cokerA} 
Let $A$ be an m-regular involutive quantale.
Then, for every surjective module homomorphism $p:P\to M$ in 
$\hbox{\rm mreg}-HMOD_A$, ${p}$ is the coker of a pair $(u, v)$ of
compact arrows from the following diagram.
\begin{equation}\label{star1}
\begin{diagram}
A^{J}&\pile{\rTo^{u}\\
\rTo_{v}}&P&\rTo^{{p}}&M.\\
\end{diagram}
\end{equation}
\end{proposition}

\begin{proof}
Note that we know e.g.\ from \cite{borceuxvi} or \cite{joyal} that in the
category $MOD_A$ of modules over $A$ the diagram 
$$
\begin{diagram}
D&\pile{\rTo^{u'}\\
\rTo_{v'}}&P&\rTo^{{p}}&M.\\
\end{diagram}
$$
exists; 
here 
$D=\{(x, y)\in P\times P: {p}(x) = {p}(y)\}$ 
is a right $A$-module that is right separated,
$u', v'$ are the respective projections.
Then ${p}$ is the coker of the pair $(u', v')$ in this category. 
Let us form the following diagram.
$$
\begin{diagram}
\hbox{\rm ess}(D)&\rTo^{\hbox{in}_D}&
D&\pile{\rTo^{u'}\\
\rTo_{v'}}&P&\rTo^{{p}}&M\\
&&&&&\rdTo_{f}&\dDashto_{g}\\
&&&&&&Z\\
\end{diagram}
$$
Here $Z$ is in $\hbox{\rm mreg}-HMOD_A$ and $f$ is a module map from $P$
to $Z$ and $g$ is defined by $g(m)=f(x)$ for $m={p}(x)$.
Let us show that our definition is correct. 
Let $m\in M$, $m={p}(x)={p}(y)$.
Then $f(x)=f(y)$ iff ${f(x)}\diamond a={f(y)}\diamond a$ for all $a\in A$ 
i.e.\ $f(x\diamond a)=f(y\diamond a)$ for all $a\in A$.
But the last condition evidently holds since $(x, y)\in D$ gives us
$(x\diamond a, y\diamond a)\in \hbox{ess}(D)$ 
i.e.\ 
$f(x\diamond a)=f(y\diamond a)$.
Since $\hbox{\rm ess}(D)$ is m-regular it is a surjective image of
$A^{\hbox{\rm\small ess}(D)}$ by the module map 
${\bf p}_{\hbox{\rm\small ess}(D)}$. 
Then we have the diagram
$$
\begin{diagram}
A^{\hbox{\rm\small ess}(D)}&
\rTo^{p_{\hbox{\rm\small ess}(D)}}&
\hbox{\rm ess}(D)&\rTo^{\hbox{in}_D}&
D&\pile{\rTo^{u'}\\
\rTo_{v'}}&P&\rTo^{{p}}&M.\\
\end{diagram}
$$
Evidently, ${p}$ is the coker of the pair 
$$
(u, v) =
(u'\circ in_D\circ p_{\hbox{\rm\small ess}(D)}, 
v'\circ in_D\circ p_{\hbox{\rm\small ess}(D)}
$$
in the category of m-regular Hilbert $A$-modules.
Recall that 
$$
u 
= 
\bigvee_{(d_1, d_2)\in \hbox{\rm\tiny ess}(D)}
u'\circ in_D\circ ((d_1, d_2)\diamond (-))
= 
\Omega_{(d_1)_{(d_1, d_2)\in \hbox{\rm\tiny ess}(D)}}
$$ 
and similarly 
$$v = \Omega_{(d_2)_{(d_1, d_2)\in \hbox{\rm\tiny ess}(D)}}$$ 
i.e.\ by Lemma \ref{cahrc} $u$ and $v$ are compact.
\end{proof}

\begin{lemma}\label{freeobje}
Let $A$ be an unital involutive quantale. 
Then
$\hbox{\rm mreg}-HMOD_{A}$ and
$\hbox{\rm mreg}-Hilb_A$ have free objects.
\end{lemma}

\begin{proof}
Let $X$ be an arbitrary set.
We put $F_A(X)=A^{X}$.
Evidently, $F_A(X)$ is an m-regular right Hilbert $A$-module. 
Let us show\ that it is free over $X$.
Evidently, the set $\{{\bf e}_x: x\in X\}$ generates $F_A(X)$ as a
submodule and we have an inclusion $\iota_X:X\to F_A(X)$ defined by
$x\mapsto {\bf e}_x$. 
By standard considerations we can check that, for any map $f:X\to M$, $M$
being an m-regular right Hilbert $A$-module, there is a unique $A$-module
map $g:F_A(X)\to M$ such that $f=g\circ i_X$.
Moreover, let us define a map $h:M\to F_A(X)$ by 
$h(n)=(\hs{f(x)}{n})_{x\in X}$.
Then, for all $(a_x)_{x\in X}\in A^{X}$ and for all $n\in M$, we have
$$
\begin{array}{lll}
\hs{g(a_x)_{x\in X}}{n}&
=&
\hs{\bigvee_{x\in X}f(x)\diamond a_x}{n}\\
&
=&
\bigvee_{x\in X}\hs{f(x)\diamond a_x}{n}\\
&
=&
\bigvee_{x\in X}{a_x}^{*}\hs{f(x)}{n}\\
&
=&
\hs{(a_x)_{x\in X}}{h(n)}.
\end{array}
$$
Then $h=g^{*}$ i.e.\ $\hbox{mreg}-Hilb_A$ has free objects.
Note that $g=\Omega_{(f(x))_{x\in X}}$.
\end{proof}

We also have the following.

\begin{lemma}\label{epion}
Let $A$ be an m-regular involutive quantale.
Assume that $M, N$ are m-regular right Hilbert $A$-modules and that
$f:M\to N$ is an epimorphism in 
${\bf C}_A \in \{\hbox{\rm mreg}-HMOD_A, \hbox{\rm mreg}-Hilb_A\}$.
Then $f(M)$ is a right Hilbert $A$-module, $f(M)\in {\bf C}_A$ and the
induced surjective map $\overline{f}:M\to f(M)$ is adjointable whenever
$f$ is. 
Moreover, we have an inner-product preserving module embedding $\hat{f}$
from $f(M)$ to $N$ such that $f(M)$ separates elements of $N$.
\end{lemma}

\begin{proof} 
Evidently, $f(M)$ is a pre-Hilbert $A$-module and, whenever $M$ is
essential then also $f(M)$ is essential.

Let us show that $f(M)$ separates $N$.
Let $n_1, n_2\in N$, $n_1\not=n_2$. 
Then
$\hs{n_1}{-}\not =\hs{n_2}{-}$
i.e.\ 
$\hs{n_1}{-}\circ f\not =\hs{n_2}{-}\circ f$ 
i.e.\
there is an element $p\in M$ such that 
$\hs{n_1}{f(p})\not =\hs{n_2}{f(p)}$. 
In particular,
$$
\hs{f(m_1)}{-}_{f(M)}=\hs{f(m_2)}{-}_{f(M)}\
\mbox{implies}\ 
f(m_1)=f(m_2)
$$
i.e.\ $f(M)$ is an m-regular right Hilbert $A$-module.
Clearly, for an adjointable $f$, $\overline{f}$ is adjointable since
$$
\hs{\overline{f}(m_1)}{f(m_2)} = 
\hs{{f}(m_1)}{f(m_2)} = 
\hs{m_1}{f^{*}(f(m_2))}.
$$

Note that the inclusion $\hat{f}:f(M) \to N$, $f(m)\mapsto f(m)$ is
evidently a module map preserving inner-product.
\end{proof}

\begin{corollary}\label{corepion}
Let $A$ be an m-regular involutive quantale.
Assume that $M, N$ are Hilbert $A$-modules and $f:M\to N$ is
an epimorphism in
${\bf C}_A \in \{\hbox{\rm mreg}-HMOD_A,
\hbox{\rm mreg}-Hilb_A\}$.
Then $f$ is a surjective map iff $\hat{f}$ is adjointable.
\end{corollary}

\begin{proof}
Evidently, if $f$ is surjective then $\hat{f}=\hbox{\rm id}_N$. 
Conversely, let $\hat{f}$ be adjointable and assume that 
$n\not\in f(M)$. 
Then, for any $m\in M$ there is an element $u_m\in M$ such that
$$
\hs{f(m)}{f(u_m)} \not= 
\hs{n}{f(u_m)}=\hs{n}{\hat{f}(f(u_m))} =
\hs{\hat{f}^{*}(n)}{f(u_m)}.$$
But $\hat{f}^{*}(n)=f(m)$ for some $m\in M$, a contradiction.
\end{proof}

\begin{corollary}\label{xcorepion}
Let $A$ be an m-regular involutive quantale.
Assume that $M, N$ are m-regular right Hilbert $A$-modules and $f$ is
an epimorphism in 
${\bf C}_A \in \{\hbox{\rm mreg}-HMOD_A, \hbox{\rm mreg}-Hilb_A\}$.
Then $f$ is a surjective map iff $f$ is a coker in ${\bf C}_A$.
\end{corollary}

\begin{proof}
Evidently, if $f$ is a surjective map so it is a coker of maps $u$ and $v$
by Proposition \ref{cokerA}.
Conversely, assume that $f$ is a coker of $u$ and $v$.
Then we have the commutative diagram.
$$
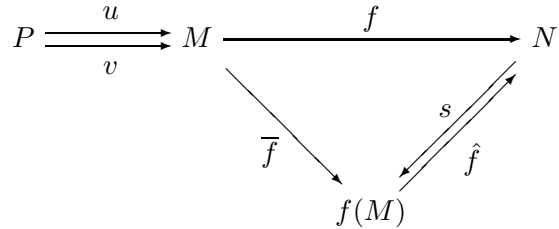
\begin{diagram}
P&\pile{\rTo^{u}\\
\rTo_{v}}&M& &\rTo^{f}& &N\\
&& &\rdTo_{\overline{f}}& &\pile{\ldDashto^{s}_{\hbox{\phantom{x}}}
\\ \ruTo_{\hat{f}}}& \\
&& & &f(M) & & \\
\end{diagram}
$$
Here $s:N\to f(M)$ is the unique module (adjointable) map such that 
$\overline{f}=s\circ f$. 
Then 
$$\hat{f}\circ (s\circ f)=\hat{f}\circ \overline{f}=f$$ 
i.e.\
$\hat{f}\circ s=\hbox{\rm id}_N$.
Hence $\hat{f}$ is onto i.e.\ $f(M)=N$.
\end{proof}

\begin{definition}\label{comma}
Let $A$ be an m-regular involutive quantale, $M, N\in{\bf C}_A$,
${\bf C}_A \in \{\hbox{\rm mreg}-HMOD_A,
\hbox{\rm mreg}-Hilb_A\}$. Assume that
$f:M\to N$ is a morphism in ${\bf C}_A$. We shall
denote by ${}_{f}{\bf C}_A$ the category with objects
all triples $\langle e, X, m\rangle$, $e:M\to X$ an
epimorphism in ${\bf C}_A$,
$m:N\to X$ an injective module map such that
$e=m\circ f$, and as arrows
$\begin{diagram}
\langle e, X, m\rangle&\rTo^{g}&\langle e', X', m'\rangle
\end{diagram}$ 
all module maps $g:X\to X'$
such that $e'=g\circ e$ and
$m'=g\circ m$ as maps. In pictures,
$$\hbox{objects}\ \langle e, X, m\rangle:
\quad\quad\quad\quad\quad\quad\quad
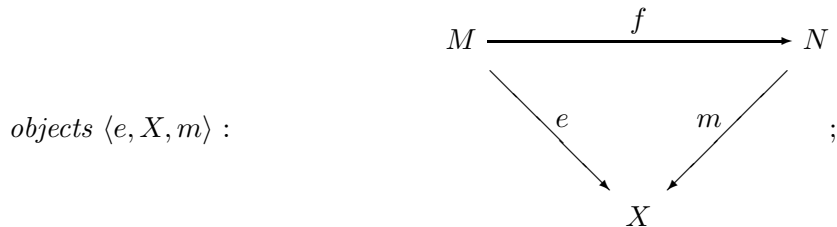
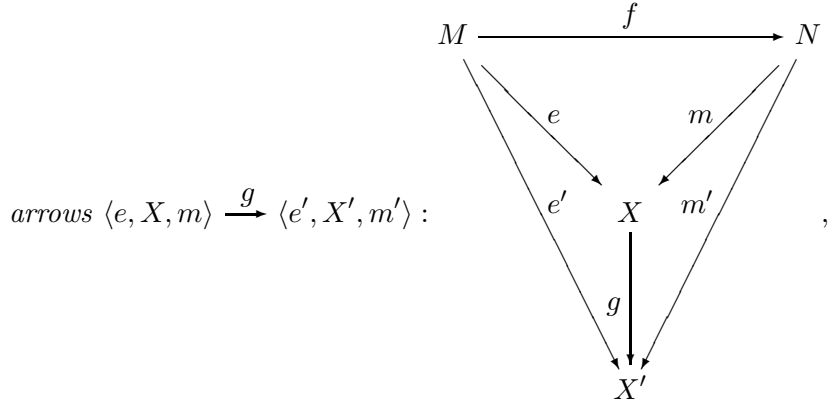
\begin{diagram}
M& &\rTo^{f}& &N\\
&\rdTo^{e}& &\ldTo^{m}& \\
& &X & & \\
\end{diagram};
$$
$$
\hbox{arrows}\ 
\begin{diagram}
\langle e, X, m\rangle&\rTo^{g}&\langle e', X', m'\rangle
\end{diagram}:
\begin{diagram}
M& &\rTo^{f}& &N\\
&\rdTo^{e}\rdTo(2,4)^{e'}& &\ldTo^{m}\ldTo(2,4)^{m'}& \\
& &X & & \\
&&\dTo^{g}&& \\
& &X' & & \\
\end{diagram},
$$
with the diagram commutative.
\end{definition}

\begin{proposition}\label{crucial}
Let $A$ be an m-regular involutive quantale, $Y, Z\in{\bf C}_A$,
${\bf C}_A \in \{\hbox{\rm mreg}-HMOD_A, \hbox{\rm mreg}-Hilb_A\}$.
Assume that $q:Y\to Z$ is an epimorphism in ${\bf C}_A$,
$\overline{q}:Y\to q(Y)$ is the induced surjective module
(adjointable) map and $\hat{q}:q(Y)\to Z$ is the induced injective
module map. 
Then $q$ is surjective if and only if the triple 
$\langle q, Z, \hat{q}\rangle$ is couniversal in 
${}_{\overline{q}}{\bf C}_A$.
\end{proposition}

\begin{proof}
Let $q$ be surjective and let $\langle e, X, m\rangle$ be in
${}_{\overline{q}}{\bf C}_A$. 
Then from the first diagram we have that $\langle q, Z, \hat{q}\rangle$ is
couniversal since $q=\overline{q}$ and $\hat{q}=\hbox{\rm id}_Z$.
Conversely, suppose that $\langle q, Z, \hat{q}\rangle$ is couniversal.
From the second diagram we have that
$\langle \overline{q}, q(Y), \hbox{\rm id}_{q(Y)}\rangle$
is couniversal as well. 
Since there is an isomorphism $h:Z\to q(Y)$ such that 
$\overline{q}=h\circ q$ we have that $q=h^{-1}\circ \overline{q}$ is
surjective.
$$
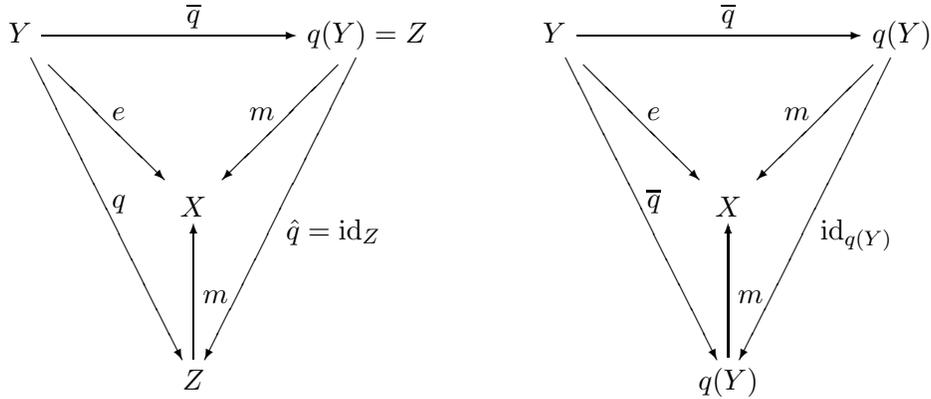
\begin{diagram}
Y& &\rTo^{\overline{q}}& &q(Y)=Z\\
&
\rdTo^{e}\rdTo(2,4)^{q}&
&
\ldTo^{m}\ldTo(2,4)_{\hat{q}=\hbox{\rm id}_Z}&
\\
& &X & & \\
&&\uTo_{m}&& \\
& &Z & & \\
\end{diagram}
\quad\quad\quad\quad
\begin{diagram}
Y& &\rTo^{\overline{q}}& &q(Y)\\
&
\rdTo^{e}\rdTo(2,4)^{\overline{q}}&
&
\ldTo^{m}\ldTo(2,4)_{\hbox{\rm id}_{q(Y)}}
&
\\
& &X & & \\
&&\uTo_{m}&& \\
& &q(Y) & & \\
\end{diagram}
$$
\end{proof}

\begin{definition}\label{defwepr} 
Let $A$ be an m-regular involutive quantale, $M\in{\bf C}_A$,
${\bf C}_A \in \{\hbox{\rm mreg}-HMOD_A, \hbox{\rm mreg}-Hilb_A\}$.
Then $M$ is called {\em weakly projective} in ${\bf C}_A$ if to every
diagram in ${\bf C}_A$
$$
\begin{diagram}
&&M\\
&&\dTo_{\varphi}\\
Y&\rTo_{q}&Z
\end{diagram}
$$
such that $q$ is a surjective morphism there is a morphism $\psi:M\to Y$
making the diagram commutative.
\end{definition}

\begin{proposition}\label{project}
Let $A$ be an unital involutive quantale.
Then $A$ is weakly projective in ${\bf C}_A$,
${\bf C}_A \in \{\hbox{\rm mreg}-HMOD_A, \hbox{\rm mreg}-Hilb_A\}$.
\end{proposition}

\begin{proof}
Assume that we have the diagram from \ref{defwepr}. 
Then there is an element $y\in Y$ such that $q(y)=\varphi(e)$. 
Define $\psi:A\to Y$ by $\psi(a)=y\diamond a$. 
Then evidently $\psi$ is an adjointable $A$-module map such that
$q\circ \psi$.
\end{proof}

\begin{corollary}\label{prproject}
Let $A$ be an unital involutive quantale.
Then $A^{J}$ is weakly projective in ${\bf C}_A$,
${\bf C}_A \in \{\hbox{\rm mreg}-HMOD_A, \hbox{\rm mreg}-Hilb_A\}$ for any
index set $J$.
\end{corollary}

\begin{proposition}
Let $A$ be an unital involutive quantale, $M\in {\bf C}_A$,
${\bf C}_A \in \{\hbox{\rm mreg}-HMOD_A, \hbox{\rm mreg}-Hilb_A\}$
Then $M$ is weakly projective in ${\bf C}_A$ if and only if $M$ is a
retract of a free Hilbert $A$-module.
\end{proposition}

\begin{proof}
Assume that $M$ is weakly projective.
Since $M$ is a surjective image of $A^{J}$ by some $q:A^{J}\to M$ we have
a morphism $\psi:M\to A^{J}$ (taking $\varphi=\hbox{\rm id}_M$) such that
$\hbox{\rm id}_M=q\circ \psi$. 
Conversely, let $M$ be a retract of a free Hilbert $A$-module $A^{J}$,
$r:A^{J}\to M$ the retraction and $i:M\to A^{J}$ the embedding such that
$r\circ i=\hbox{\rm id}_M$. 
Let $Y, Z$ be Hilbert $A$-modules, $q:Y\to Z$ a surjective morphism,
$\phi:M\to Z$ a morphism. 
Since $A^{J}$ is weakly projective there is a morphism $\psi':A^{J}\to Y$
such that $q\circ \psi'=\varphi\circ r$. 
Let us put $\psi=\psi'\circ i$. 
Then
$$
q\circ\psi=q\circ\psi'\circ i=\varphi\circ r\circ i=\varphi.
$$
\end{proof}

\section{The Eilenberg-Watts and fundamental Morita theorems for Hilbert
modules}

\begin{definition}\label{defome} 
We say that two m-regular involutive quantales $A$ and $B$ are {\em Morita
equivalent} if there exist $*$-functors 
$F : \hbox{\rm mreg}-{Hilb}_A \rightarrow \hbox{\rm mreg}-Hilb_B$
and 
$G : \hbox{\rm mreg}-{Hilb}_B \rightarrow \hbox{\rm mreg}-Hilb_A$,
such that $F G \cong Id$ and $G F \cong Id$ unitarilly.
Such $F$ and $G$ will be called {\em equivalence functors}.
\end{definition}

\begin{theorem}[Eilenberg-Watts theorem for m-regular involutive
quantales]\label{EW}
Let $A$ and $B$ be m-regular involutive quantales and let $F$ be a
colimit-preserving *-functor from
$\hbox{\rm mreg}-Hilb_A$ to $\hbox{\rm mreg}-Hilb_B$.
Then there is a right Hilbert $A-B$ bimodule $Z$ such that $F(-)$ is
naturally unitarily isomorphic to the interior tensor product 
$(-) \dot{\otimes}_{A} Z$. 
That is, there is a natural isomorphism between these functors, which
implements a unitary isomorphism $F(Y) \cong Y \dot{\otimes}_{A} Z$ 
for all $Y \in \hbox{\rm mreg}-Hilb_A$.
\end{theorem}

\begin{proof}
Define $Z = F(A)$. 
Then $Z$ is an m-regular right Hilbert $B$-module. 
We make $Z$ into a left $A$-module by defining $a\bullet z = F(L(a))(z)$,
for $a \in A$, $z \in Z$. 
Here $L(a) : A \rightarrow A$ is the adjointable map $b \mapsto ab$. 
Let us check that $Z$ is a right Hilbert $A-B$ bimodule. 
We have
$$
\begin{array}{lll}
\hs{a\bullet x}{y}&
=&
\hs{F(L(a))(x)}{y}\\
&
=&
\hs{x}{F(L(a))^{*}(y)}\\
&
=
&\hs{x}{F(L(a)^{*})(y)}\\
&
=&
\hs{x}{F(L(a^{*}))(y)}\\
&
=&
\hs{x}{a^{*}\bullet y}.
\end{array}
$$

Now, we shall check that $Z$ is essential with respect to $\bullet$. 
Note that we know that the map ${\bf p}:A^{A}\to A$,
${\bf p}=\bigvee_{a\in A}L(a)\circ \pi_a$ is surjective 
i.e.\ it is a coker of adjointable maps in $\hbox{rm mreg}-Hilb_A$. 
Then $F({\bf p}):F(A)^{A}\to F(A)$ is a coker in 
$\hbox{\rm mreg}-{Hilb}_B$ i.e.\
$$
F({\bf p}) = 
\bigvee_{a\in A}F(L(a))\circ F(\pi_a)=
\bigvee_{a\in A}a\bullet F(\pi_a)
$$ 
is surjective
by Corollary \ref{xcorepion}.

It follows that, for all $Y \in \hbox{\rm mreg}-{Hilb}_A$, the interior
tensor product $Y \dot{\otimes}_{A} Z$ is well defined. 
We define an adjointable map
$\tau_Y : Y \dot{\otimes}_{A} Z \rightarrow F(Y)$ 
by
$\tau_Y(y \dot{\otimes}_A z) = F(L(y))(z)$,
where $L(y) : A \rightarrow Y$ is the map $L(y)(a) = y\diamond a$.
Let us prove that $\tau_Y$ is a unitary natural transformation.

First note that for $Y = A^{J}$ this is easy, in fact $\tau_{A^{J}}$ is
the canonical unitary isomorphism from 
$$
A^{J} \dot{\otimes}_{A} F(A) \cong 
(A \dot{\otimes}_{A} F(A))^{J} \cong
F(A)^{J} \cong 
F(A^{J})
.
$$
This is because if $y = (a_j)_{j\in J} \in A^{J}, z\in F(A)$ then
$$
\begin{array}{lll}
\multicolumn{3}{l}{
(a_j)_{j\in J}\dot{\otimes}_{A} z \mapsto
(a_j\dot{\otimes}_{A} z)_{j\in J} \mapsto 
(a_j\bullet z)_{j\in J}
}
\\
&
=&
(F(L(a_j))(z))_{j\in J} \mapsto 
\bigvee_{j\in J} F(i_j)(F(L(a_j))(z))\\
&
=&
\bigvee_{j\in J} F(i_j \circ L(a_j))(z)\\
&
=&
\bigvee_{j\in J} F(L(i_j(a_j)))(z)\\
&
=&
F(\bigvee_{j\in J} L(i_j(a_j)))(z)\\
&
=&
F(L(\bigvee_{j\in J} i_j(a_j)))(z)\\
&
=&
F(L(y))(z)\\
&
=&
\tau_Y(y \dot{\otimes}_A z).
\end{array}
$$

Now fix $Y \in \hbox{\rm mreg}-{Hilb}_A$. 
We have then the following commutative diagram for a suitable index set 
$J$ such that ${\bf p}_Y$ is the coker of adjointable maps $u, v$.
\begin{equation}\label{star2}
\begin{diagram}
A^{J}&\pile{\rTo^{u}\\ 
\rTo_{v}}&A^{Y}&\rTo^{{\bf p}_Y}&Y\\
\end{diagram}
\end{equation}

Applying $F$, $(-)\dot{\otimes}_A Z$ and again $F$ we get a diagram.
\begin{equation}\label{star3}
\begin{diagram}
F(A^{J})&\pile{\rTo^{F(u)}\\ \rTo_{F(v)}}&F(A^{Y})&
\rTo^{F({\bf p}_Y)}&F(Y)\\
\dTo^{\tau_{A^{J}}^{*}}& &\dTo_{\tau_{A^{Y}}^{*}}&
&\dTo_{\tau_{Y}^{*}}\\
A^{J}\dot{\otimes}_A Z&
\pile{\rTo^{u\dot{\otimes}_A \hbox{\rm id}_Z}\\ 
\rTo_{v\dot{\otimes}_A \hbox{\rm id}_Z}}&
A^{Y}\dot{\otimes}_A Z&
\rTo^{{\bf p}_Y\dot{\otimes}_A \hbox{\rm id}_Z}&
Y\dot{\otimes}_A Z\\
\dTo^{\tau_{A^{J}}}& &\dTo_{\tau_{A^{Y}}}& &\dTo_{\tau_{Y}}\\
F(A^{J})&\pile{\rTo^{F(u)}\\ \rTo_{F(v)}}&F(A^{Y})&
\rTo^{F({\bf p}_Y)}&F(Y)\\
\end{diagram}
\end{equation}

We shall first check that the diagram is commutative.

Note that 
$$\begin{array}{lll}
\tau_Y({\bf p}_Y(w) \dot{\otimes}_A z)&
=&
F(L({{\bf p}_Y(w)}))(z)\\
&
=&
 F({\bf p}_Y) F(L(w)) (z)\\
&
=&
F({\bf p}_Y) \tau_{A^{Y}} (w \dot{\otimes}_A z)
\end{array}$$
for $w \in A^Y$, $z \in Z$. 
Similarly, let $y\in Y$, $w=(a_x)_{x\in Y} \in A^Y$, $z, z' \in Z$. 
Then
$$
\begin{array}{lll}
\langle
{\bf p}_Y\dot{\otimes}_A \hbox{\rm id}_Z ((a_x)_{x\in Y}\dot{\otimes}_A z'), 
{y\dot{\otimes}_A z}
\rangle
&
=&
\hs{z'}{\hs{{\bf p}_Y(a_x)_{x\in Y}}{y}\bullet z}\\
&
=&
\hs{\hs{y}{{\bf p}_Y(a_x)_{x\in Y}}\bullet z'}{z}\\
&
=&
\bigvee_{x\in Y}\hs{\hs{y}{x\diamond a_x}\bullet z'}{z}\\
&
=&
\bigvee_{x\in Y}\hs{F(L(\hs{y}{x\diamond a_x}))(z')}{z}
\end{array}
$$
and
$$
\begin{array}{lll}
\multicolumn{3}{l}{
\langle
( \tau_{Y}^{*}\circ F({\bf p}_Y)\circ \tau_{A^{Y}} )
((a_x)_{x\in Y}\dot{\otimes}_A z')
, 
{y\dot{\otimes}_A z}
\rangle
}
\\
&
=&
\hs
{(F({\bf p}_Y) \circ \tau_{A^{Y}})((a_x)_{x\in Y}\dot{\otimes}_A z')}
{\tau_Y(y\dot{\otimes}_A z)}
\\
&
=&
\hs
{
(F(L(y)^{*}) \circ F({\bf p}_Y) \circ \tau_{A^{Y}})
((a_x)_{x\in Y}\dot{\otimes}_A z')
}
{z}
\\
&
=&
\hs
{\bigvee_{x\in Y} (F(L(y)^{*}\circ {\bf p}_Y \circ i_x \circ L(a_x))(z')}
{z}
\\
&
=&
\bigvee_{x\in Y}
\hs
{(F((\hs{y}{-}) \circ (x\diamond (-)) \circ (a_x\cdot (-)))(z')}
{z}
\\
&
=&
\bigvee_{x\in Y}
\hs
{(F(L(\hs{y}{x\diamond a_x}))(z')}
{z}.
\end{array}
$$

Then, since the upper and lower horizontal lines of the diagram (\ref{star3})
produce the respective cokers $F({\bf p}_Y)$,
the composite right square of this diagram is a pushout in $\hbox{\rm 
mreg}-{Hilb}_B$.
\begin{equation}\label{star4}
\begin{diagram}
F(A^{Y})&\rTo^{F({\bf p}_Y)}&F(Y)\\
\dTo^{\hbox{\rm id}_{F(A^{Y})}}&
&\dTo_{\tau_Y\circ\tau_{Y}^{*}}\\
F(A^{Y})&\rTo^{F({\bf p}_Y)}&F(Y)\\
\end{diagram}
\end{equation}
i.e.\
$\tau_Y\circ \tau_{Y}^{*}=\hbox{\rm id}_{F(Y)}$.
Similarly, we have this commutative diagram.
\begin{equation}\label{star5}
\begin{diagram}
A^{Y}\dot{\otimes}_A Z&
\rTo^{{\bf p}_Y\dot{\otimes}_A \hbox{\rm id}_Z}&
Y\dot{\otimes}_A Z\\
\dTo_{\hbox{\rm id}_{A^{Y}\dot{\otimes}_A Z}}&
&\dTo^{\hbox{\rm id}_{Y\dot{\otimes}_A Z}}
\dTo_{\tau_Y^{*}\circ \tau_{Y}}\\
A^{Y}\dot{\otimes}_A Z&
\rTo^{{\bf p}_Y\dot{\otimes}_A \hbox{\rm id}_Z}&
Y\dot{\otimes}_A Z\\
\end{diagram}
\end{equation}

Since ${\bf p}_Y\dot{\otimes}_A \hbox{\rm id}_Z$ is a surjective 
adjointable map it is an epimorphism and we have
$\tau_Y^{*}\circ \tau_Y=\hbox{\rm id}_{Y\dot{\otimes}_A Z}$.

Let $f:X\to Y$ be a morphism in $\hbox{\rm mreg}-{Hilb}_A$, $x\in X$ and
$z\in Z$. 
Then 
$$
\begin{array}{lll}
(\tau_Y\circ (f\dot{\otimes}_A \hbox{\rm id}_Z))(x\dot{\otimes}_A z)&
=&
F(L(f(x)))(z)\\
&
=&
F(f\circ L(x))(z)\\
&
=&
F(f)(F(L(x))(z))\\
&
=&
(F(f)\circ \tau_X)(x\dot{\otimes}_A z)
\end{array}
$$ 
and since $\tau_X$, $\tau_Y$ are unitary also 
$$
\tau_Y^{*}\circ F(f) = (f\dot{\otimes}_A \hbox{\rm id}_Z))\tau_X^{*}
$$ 
i.e.\ $\tau_Y, \tau_Y^{*}$ are natural transformations.
\end{proof}

\begin{lemma}\label{bassth}
Let $A$ and $B$ be m-regular involutive quantales and let $X_1$ and $X_2$
be essential right Hilbert $A-B$ bimodules such that
$F_1=(-)\dot{\otimes}_A X_1$ and
$F_2=(-)\dot{\otimes}_A X_2$ 
are functors from
$\hbox{\rm mreg}-{Hilb}_A$ to
$\hbox{\rm mreg}-{Hilb}_B$. 
Then $F_1$ and $F_2$ are equivalent if and only if $X_1\cong X_2$
unitarily and as bimodules.
\end{lemma}

\begin{proof}
Assume that $\alpha:F_1\to F_2$
is a natural isomorphism such that $\alpha_{M}:F_1(M)\to F_2(M)$
is a unitary map. Then we have the following diagram.
$$
\begin{diagram}
A\dot{\otimes}_A X_1&\rTo^{\alpha_{A}}&A\dot{\otimes}_A X_1\\
\uTo^{\kappa_{X_1}}& &\dTo^{\chi_{X_2}}\\
X_1&\rTo^{h}&X_2\\
\end{diagram}
$$
Here $\chi_{X_2}(a\dot{\otimes}_{A} x)=a\bullet x$,
$\kappa_{X_1}=\chi_{X_1}^{*}$ and
$h=\chi_{X_2}\circ \alpha_{A}\circ \kappa_{X_1}$.
Left multiplication in $A$ are adjointable
maps, so a natural transformation $\alpha_A$
has to preserve them. From \cite{pasekainterior} we know that
$\kappa_{X_1}$ and $\chi_{X_2}$ are unitary maps and
bimodule homomorphisms.
Hence $h$ is unitary and bimodule homomorphism. The other
direction is evident.
\end{proof}

The following theorem is an involutive quantale version of Morita's 
fundamental theorem.

\begin{theorem}\label{xxmainmorita} 
Let $A$ and $B$ be m-regular involutive quantales. 
Then $A$ and $B$ are Morita-equivalent if and only if there are an 
essential right Hilbert $A-B$ bimodule $X$ and an essential right Hilbert
$B-A$ bimodule $Y$ such that $X\dot{\otimes}_B Y\cong A$ and 
$Y\dot{\otimes}_A X\cong B$ as bimodules.
\end{theorem}

\begin{proof}
We have $F(-) \cong - \dot{\otimes}_{A} X$
and $G(-) \cong - \dot{\otimes}_{B} Y$.
Composing these two we obtain 
$$
A \cong 
GF(A) \cong 
A \dot{\otimes}_{A} X \dot{\otimes}_{B} Y \cong 
X \dot{\otimes}_{B} Y.
$$ 
Similarly,
$Y \dot{\otimes}_{A} X \cong B$,
and this identification and the last are unitarily, and as bimodules, the
latter exactly as in pure algebra \cite{bass}.

Conversely, let $\lambda:X\dot{\otimes}_B Y\to A$ be the unitary 
isomorphism, 
$F(-) \cong - \dot{\otimes}_{A} X$
and
$G(-) \cong - \dot{\otimes}_{B} Y$.
Then
$$(GF)(M)\cong M\dot{\otimes}_{A} X \dot{\otimes}_{B} Y.$$
Let us define
$$\tau_M=\pi_M\circ (\hbox{\rm id}_M\dot{\otimes}_A \lambda),$$
where $\pi_M:M\dot{\otimes_A} A\to M$
is the canonical unitary isomorphism defined by
$m\dot{\otimes_A} a\mapsto m\diamond a$. 
We get a natural isomorphism 
$\tau:GF\to \hbox{\rm Id}_{\hbox{\tiny\rm mreg}-Hilb_A}$. 
By symmetry, $A$ is Morita equivalent to $B$.
\end{proof}

\begin{corollary}\label{projmor}
Let $A$ and $B$ be unital involutive quantales. 
Then $A$ and $B$ are Morita-equivalent if and only if their categories of
weakly projective m-regular right Hilbert modules are equivalent.
\end{corollary}

\begin{proof} 
It follows from the fact that weak projectivity is a categorical property. 
Hence weakly projective objects are mapped on weakly projective objects. 
Conversely, any equivalence between categories of weakly projective 
m-regular right Hilbert modules can be easily extended to an equivalence
of m-regular right Hilbert modules.
\end{proof}

We recall that if $A$ and $B$ are m-regular involutive quantales then an
imprimitivity Hilbert $A-B$ bimodule is an $A-B$ bimodule $X$, which is a
full m-regular right Hilbert $B$-module, and also a full m-regular left
Hilbert $A$-module, such that
$\; _A\langle x, y \rangle z = x \langle y, z \rangle_B$ 
whenever $x,y,z \in X$. 
The existence of such an $X$ is the definition of $A$ and $B$ being {\em
strongly Morita equivalent}. 
Our theorem gives a functorial characterization of such an $X$, and of
strong Morita equivalence of $A$ and $B$.

\begin{theorem}\label{mainmorita} 
Let $A$ and $B$ be m-regular involutive quantales.
Suppose that 
$F\colon \hbox{\rm mreg}-{Hilb}_A\rightarrow \hbox{\rm mreg}-{Hilb}_B$,
$G\colon \hbox{\rm mreg}-{Hilb}_B\rightarrow \hbox{\rm mreg}-{Hilb}_A$
are equivalence *-functors, with $FG \cong Id$ and $GF \cong Id$ via
unitary natural isomorphisms.
Then $A$ and $B$ are strongly Morita equivalent. 
Moreover, the $A-B$ bimodule $X = F(A)$ from the Eilenberg-Watts theorem
is an imprimitivity bimodule implementing the strong Morita
equivalence. 
As in that theorem, $F(-) \cong - \dot{\otimes}_{A} X$ naturally and
unitarily.
Conversely, any $A-B$ imprimitivity bimodule $X$ implements such a
functorial isomorphism.
\end{theorem}

\begin{proof}
It follows from the fact that $X$ is an imprimitivity $A-B$ bimodule
if and only if there is an essential right Hilbert $B-A$ bimodule $Y$
such that $X\dot{\otimes}_B Y\cong A$ and $Y\dot{\otimes}_A X\cong {}B$ as
right-Hilbert bimodules (see \cite{pasekainterior}) and Theorem
\ref{xxmainmorita}.
\end{proof}

\begin{proposition}\label{isoclass}
Let $A$, $B$ and $C$ be m-regular involutive quantales.
Then the unitary isomorphism classes of equivalence functors
$\hbox{\rm mreg}-{Hilb}_A \rightarrow \hbox{\rm mreg}-{Hilb}_B$
are in a 1-1 correspondence with the unitary equivalence classes of
imprimitivity $A-B$ bimodules. 
Composition of such functors corresponds to the interior tensor product of
the bimodules.
\end{proposition}

\begin{proof} 
Every imprimitivity $A-B$ bimodule $X$ gives rise to an equivalence
$(-)\dot{\otimes}_{A} X: 
\hbox{\rm mreg}-{Hilb}_A\to \hbox{\rm mreg}-{Hilb}_B$, 
the isomorphism type depends only on the isomorphism type of $X$ (as we 
saw in Lemma \ref{bassth}). 
Conversely, if
$F:\hbox{\rm mreg}-{Hilb}_A\to \hbox{\rm mreg}-{Hilb}_B$ 
is an equivalence, $F(A)$ is an imprimitivity $A-B$ bimodule and its
isomorphism type depends only on that of $F$. 
If $Y$ is an imprimitivity $B-C$ bimodule of an equivalence functor, the
composition of the equivalences
$\hbox{\rm mreg}-{Hilb}_A \to 
\hbox{\rm mreg}-{Hilb}_B \to 
\hbox{\rm mreg}-{Hilb}_C$ 
is given by
$(-)\dot{\otimes}_{A} (X\dot{\otimes}_B Y)$. 
In particular, $X\dot{\otimes}_B Y$ is an imprimitivity $A-C$ 
bimodule.
\end{proof}

\begin{corollary}\label{selfeq}
Let $A$ be an m-regular involutive quantale.
Then the unitary isomorphism classes of self-equivalence functors
$\hbox{\rm mreg}-{Hilb}_A \rightarrow \hbox{\rm mreg}-{Hilb}_A$
under a composition form a group isomorphic to the group of the unitary
equivalence classes of imprimitivity $A-A$ bimodules.
\end{corollary}

\begin{corollary}\label{dual}
Let $A$ and $B$ be m-regular involutive quantales. 
Then $A$ and $B$ are Morita-equivalent if and only if $A^{d}$ and $B^{d}$
are Morita-equivalent.
\end{corollary}

\begin{proof} 
Let $A$ and $B$ are Morita-equivalent. 
Then we have an imprimitivity $A-B$ bimodule $X$. 
Then $Y$ with the same order and inner products as $X$ and equipped with
operations $a\bullet_d x=a^{*}\bullet x$ and 
$x\diamond_{d} b=x\diamond x^{*}$ is
evidently an imprimitivity $A^{d}-B^{d}$ bimodule.
\end{proof}

\begin{proposition}\label{matreq} 
Let $A$ be an m-regular involutive quantale. 
Then $A$ is Morita equivalent to the matrix quantale $\mathcal{M}^{J}(A)$.
\end{proposition}

\begin{proof} 
Note that $\mathcal{M}^{J}(A)\cong \mathcal{K}_{A}(A^{J}, A^{J})$.
\end{proof}

\begin{proposition}\label{xmatreq} 
Let $S$ be a sup-semilatticce with a duality. 
Then the involutive subquantale $\mathcal{Q}_{0}(S)$ of $\mathcal{Q}(S)$
that is generated by right-sided elements of $\mathcal{Q}(S)$ is Morita
equivalent to the $2$-element Boolean algebra ${\bf 2}$.
In particular, $\mathcal{M}^{J}({\bf 2})$ is Morita equivalent to ${\bf 2}$.
\end{proposition}

\begin{proof}
Recall that $\mathcal{Q}_{0}(S)\cong \mathcal{K}_{{\bf 2}}(S, S)$.
\end{proof}

\begin{lemma}\label{expl}
Let $E$ be an imprimitivity Hilbert $B$-$A$ bimodule with 
$A=\mathcal{K}_{\bf 2}(S_A)$, $S_A$ being a sup-lattice with a duality. 
Then there is a sup-lattice with a duality $S_B$ such that 
$B\cong \mathcal{K}_{\bf 2}(S_B)$ and $E\cong \mathcal{K}_{\bf 2}(S_A, S_B)$.
\end{lemma}

\begin{proof} 
We write $\iota$ for the identity representation of $\mathcal{K}_{\bf
2}(S_A)$ on $S_A$ and put $S_B=E\dot{\otimes}_{A} S_A$. 
Left multiplication in the first component of $S_B$ then defines a
faithful representation of $B\simeq \mathcal{K}_{A}(E)$ onto $S_B$ 
i.e.\
$B\simeq \mathcal{K}_{\bf 2}(S_B)$
(see \cite[Corollary 1.13]{pasekainterior}).

Now we shall define a mapping
$$
\Theta \colon E\to \mathcal{A}_{\bf 2}(S_A , S_B),\
x \mapsto (s \mapsto x \dot{\otimes}_{A} s ).
$$
Namely, 
$$
\hs{\Theta(x)(s)}{y\dot{\otimes}_{A} t}
=
\hs{x\dot{\otimes}_{A} s}{y \dot{\otimes}_{A} t}
=
\hs{s}{\hs{x}{y}_{A}\bullet t}
$$
i.e.\
$$
\Theta(x)^{*}(y\dot{\otimes}_{A} t) = \hs{x}{y}_{A}\bullet t.
$$
We also observe that
$$
\Theta (x)^{*}\Theta (y)=\rip{A}{x}{y}.
$$ 
If $a\in A$, then
$$
\Theta (x\diamond a)(s) =(x\diamond a)\dot{\otimes}_{A} s 
=
x \dot{\otimes}_{A} a\bullet s
,
$$ 
so that 
$$
\Theta(E) 
= 
\Theta (E\diamond A)
\subseteq 
\Theta(E)\circ \mathcal{K}_{\bf 2}(S_A)
\subseteq 
\mathcal{K}_{\bf 2}(S_A, S_B)
.
$$ 
Similarly,
$$
\Theta_{y\dot{\otimes}_A t, s}
=
y\dot{\otimes}_A t\diamond \hs{s}{-}
=
\Theta(y)\circ \Theta_{t, s}
$$
i.e.\
$$
\mathcal{K}_{\bf 2}(S_A, S_B)
\subseteq 
\Theta(E)\circ \mathcal{K}_{\bf 2}(S_A)
\subseteq 
\Theta(E)
.
$$
So we conclude that
$$
\Theta (E)
=
\mathcal{K}_{\bf 2}(S_A, S_B)
.
$$ 
It is a straightforward computation that $\Theta$ is an {isomorphism} of
Hilbert bimodules.
\end{proof}

\begin{lemma}\label{misa} 
For any m-regular involutive quantale $A$,
$$
\mathcal{A}_A(A)\cap {}_{A}\mathcal{A}(A)= HMOD_A(A)\cap {}_A{}HMOD(A).
$$
\end{lemma}

\begin{proof} 
It is enough to show that any bimodule endomorphism $f:A\to A$ is
adjointable. 
Note that
$$
\begin{array}{lll}
\hs{m}{f(n)}&
=&
m^{*}\cdot f(n)\\
&
=&
f(m^{*}\cdot n)\\
&
=&
f(m^{*})\cdot n\\
&
=&
(f(m^{*})^{*})^{*}\cdot n\\
&
=&
\hs{f(m^{*})^{*}}{n}\mu\ \mbox{for all $m$, $n\in A$}
\end{array}
$$ 
i.e.\ $f$ has an adjoint.
\end{proof}

\begin{definition}\label{natgro}
Let $A$ be an m-regular involutive quantale. 
The set of adjointable natural transformations from the identity functor
$\hbox{\rm Id}_{\hbox{\small\rm mreg}-Hilb_A}$ to itself is called
$\hbox{\rm ANat}(A)$. 
It is a unital commutative involutive quantale with the composition of
natural transformations as multiplication, with the adjoint of a natural
transformation as involution
$$
\begin{diagram}
M&\rTo^{\eta_{M}}&M\\
\dTo^{f}&&\dTo^{f}\\
N&\rTo^{\eta_{N}}&N\\
\end{diagram}
\longmapsto
\begin{diagram}
M&\rTo^{\eta_{M}^{*}}&M\\
\uTo^{f^{*}}&&\uTo^{f^{*}}\\
N&\rTo^{\eta_{N}^{*}}&N\\
\end{diagram}
$$
and with the join given by
$$
\begin{diagram}
M&\rTo^{{\eta_{i}}_{M}}&M\\
\dTo^{f}&&\dTo^{f}\\
N&\rTo^{{\eta_{i}}_{N}}&N\\
\end{diagram}, 
i\in I\quad
\longmapsto\quad
\begin{diagram}
M&\rTo^{\bigvee_{i}{\eta_{i}}_{M}}&M\\
\dTo^{f}&&\dTo^{f}\\
N&\rTo^{\bigvee_{i}{\eta_{i}}_{N}}&N\\
\end{diagram}.
$$
The {\em centre} of $A$ is the unital commutative involutive quantale
$\hbox{\rm Cen}(A)=\mathcal{A}_{A}(A)\cap {}_{A}\mathcal{A}(A)$.
\end{definition}

The elements of $\hbox{\rm Cen}(A)$ are by Lemma \ref{misa} $A$-bimodule
endomorphisms on $A$.
Note that evidently $\hbox{\rm Cen}(A)$ is commutative since, for all
$a\in A$, $a=\bigvee_{j}a_j\cdot b_j$, $f$, $g\in \hbox{\rm Cen}(A)$, we
have
$$
\begin{array}{lll}
(f\circ g)(a)&
=&
(f\circ g)(\bigvee_{j}a_j\cdot b_j)\\
&
=&
\bigvee_{j}f(g(a_j\cdot b_j))\\
&
=&
\bigvee_{j}f(g(a_j)\cdot b_j)\\
&
=&
\bigvee_{j}g(a_j)\cdot f(b_j)\\
&
=&
\bigvee_{j}g(f(a_j\cdot b_j))\\
&
=&
(g\circ f)(\bigvee_{j}a_j\cdot b_j)\\
&
=&
(g\circ f)(a).
\end{array}
$$
Whenever $A$ is commutative we can see $A$ as a complete ${}^{*}$-ideal of
$\hbox{\rm Cen}(A)$ by identifying an element $a\in A$ with the bimodule
endomorphism on $A$ induced by multiplication by $a$.

The following theorem is based on the theorem 4.2 from \cite{ara} for
involutive rings.

\begin{theorem}\label{aracoriimp}
Let $A$, $B$ be Morita equivalent m-regular involutive quantales. 
Then
\begin{enumerate}
\item 
$\hbox{Cen}(A)$ and $\hbox{Cen}(B)$ are isomorphic as involutive
quantales.
\item 
If $A$ and $B$ are commutative then $A$ and $B$ are isomorphic as
involutive quantales.
\end{enumerate}
\end{theorem}

\begin{proof}
(1). 
By Theorem \ref{mainmorita} there is an $A-B$ imprimitivity bimodule $X$
implementing the Morita equivalence. 
We shall now define a map $\gamma:\hbox{Cen}(A)\to \hbox{Cen}(B)$ by the
prescription
$$
\gamma(f)(\bigvee_{i}\rip{B}{x_i}{r_i\bullet y_i}) =
\bigvee_{i}\rip{B}{x_i}{f(r_i)\bullet y_i}
$$
for all $x_i, y_i\in X$ and $r_i\in A$. 
Let us check that $\gamma$ is well defined. 
Assume that 
$$
\bigvee_{i}\rip{B}{x_i}{r_i\bullet y_i} =
\bigvee_{j}\rip{B}{u_j}{p_j\bullet v_j}.
$$ 
Let $z, t\in X$ and $r'\in A$. 
Then
$$
\begin{array}{lll}
\rip{B}{z}{r'\bullet t} \cdot 
(\bigvee_{i}\rip{B}{x_i}{f(r_i)\bullet y_i})&
=&
\bigvee_{i}
\rip{B}{z}{(r'\cdot\rip{A}{t}{x_i}\cdot f(r_i)) \bullet y_i}
\\
&
=&
\bigvee_{i}
\rip{B}{z}{(f(r') \cdot \rip{A}{t}{x_i}\cdot r_i) \bullet y_i}
\\
&
=&
\rip{B}{z}{f(r')\bullet t} \cdot 
(\bigvee_{i}\rip{B}{x_i}{r_i\bullet y_i})
\\
&
=&
\rip{B}{z}{f(r')\bullet t} \cdot 
(\bigvee_{j}\rip{B}{u_j}{p_j\bullet v_j})
\\
&
=
&
\rip{B}{z}{r'\bullet t} \cdot 
(\bigvee_{j}\rip{B}{u_j}{f(p_j)\bullet v_j}).
\end{array}
$$

We show that $\gamma(f)\in \hbox{\rm Cen}(B)$. 
Let $b_j\in B$, 
$$b_j=\bigvee_{i\in I_j}\rip{B}{x^{j}_i}{r^{j}_i\bullet y^{j}_i}.$$
Then
$$
\begin{array}{lll}
\gamma(f)(\bigvee_{j} b_j)&
=&
\gamma(f)(\bigvee_{j}
\bigvee_{i\in I_j}\rip{B}{x^{j}_i}{r^{j}_i\bullet y^{j}_i})
\\
&
=&
\bigvee_{j}
\bigvee_{i\in I_j}\rip{B}{x^{j}_i}{f(r^{j}_i)\bullet y^{j}_i}
\\
&
=&
\bigvee_{j}
\gamma(f) (\bigvee_{i\in I_j}\rip{B}{x^{j}_i}{{j}_i\bullet y^{j}_i})
\\
&
=&
\bigvee_{j}\gamma(f)(b_j)
\end{array}
$$
i.e.\ $\gamma(f)$ is a sup-lattice homomorphism.
Now, let $b$, $c\in B$, 
$$b=\bigvee_{i}\rip{B}{x_i}{r_i\bullet y_i}.$$
Then
$$
\begin{array}{lll}
\gamma(f)(b\cdot c)&
=&
\gamma(f)(\bigvee_{i}\rip{B}{x_i}{r_i\bullet y_i\diamond c})\\
&
=&
\bigvee_{i}\rip{B}{x_i}{f(r_i)\bullet y_i\diamond c}\\
&
=&
\bigvee_{i}\rip{B}{x_i}{f(r_i)\bullet y_i}\cdot c\\
&
=&
\gamma(f)(b)\cdot c
\end{array}
$$
and
$$
\begin{array}{lll}
\gamma(f)(c\cdot b)&
=&
\gamma(f)(\bigvee_{i}\rip{B}{x_i\diamond c^{*}}{r_i\bullet y_i})\\
&
=&
\bigvee_{i}\rip{B}{x_i\diamond c^{*}}{f(r_i)\bullet y_i}\\
&
=&
c\cdot \bigvee_{i}\rip{B}{x_i}{f(r_i)\bullet y_i}\\
&
=&
c\cdot \gamma(f)(b)
\end{array}
$$
i.e.\ $\gamma(f)$ is an bimodule endomorphism.

Let us check that $\gamma$ is an involutive quantale homomorphism from
$\hbox{\rm Cen}(A)$ to $\hbox{\rm Cen}(B)$. 
Let $f, g, f_j\in \hbox{\rm Cen}(A)$, $j\in J$. 
Then
$$
\begin{array}{lll}
\gamma(\bigvee_{j}f_j)
(\bigvee_{i}\rip{B}{x_i}{r_i\bullet y_i})&
=&
\bigvee_{i}\rip{B}{x_i}{\bigvee_{j}f_j(r_i)\bullet y_i}\\
&
=&
\bigvee_{j}\bigvee_{i}\rip{B}{x_i}{f_j(r_i)\bullet y_i}\\
&
=&
\bigvee_{j}\gamma(f_j)
(\bigvee_{i}\rip{B}{x_i}{r_i\bullet y_i}),
\end{array}
$$
$$
\begin{array}{lll}
\gamma(f^{*})
(\bigvee_{i}\rip{B}{x_i}{r_i\bullet y_i})&
=&
\bigvee_{i}\rip{B}{x_i}{f^{*}(r_i)\bullet y_i}\\
&
=&
\bigvee_{i}\rip{B}{x_i}{f(r_i^{*})^{*}\bullet y_i}\\
&
=&
\bigvee_{i}\rip{B}{f(r_i^{*})\bullet x_i}{y_i}\\
&
=&
(\bigvee_{i}\rip{B}{y_i}{f(r_i^{*})\bullet x_i})^{*}\\
&
=&
(\gamma(f)(\bigvee_{i}\rip{B}{y_i}{r_i^{*}\bullet x_i}))^{*}\\
&
=&
\gamma(f)^{*}(\bigvee_{i}\rip{B}{y_i}{r_i^{*}\bullet x_i}^{*})\\
&
=&
\gamma(f)^{*}(\bigvee_{i}\rip{B}{r_i^{*}\bullet x_i}{y_i})\\
&
=&
\gamma(f)^{*}(\bigvee_{i}\rip{B}{x_i}{r_i\bullet y_i})
\end{array}
$$
and
$$
\begin{array}{lll}
\gamma(f\circ g)
(\bigvee_{i}\rip{B}{x_i}{r_i\bullet y_i})&
=&
\bigvee_{i}\rip{B}{x_i}{(f\circ g)(r_i)\bullet y_i}\\
&
=&
\bigvee_{i}\rip{B}{x_i}{f(g(r_i)\bullet y_i}\\
&
=&
\gamma(f)(\bigvee_{i}\rip{B}{x_i}{g(r_i)\bullet y_i})\\
&
=&
(\gamma(f)\circ \gamma(g))
(\bigvee_{i}\rip{B}{x_i}{\bigvee_{j}r_i\bullet y_i}).
\end{array}
$$

Similarly, we have an involutive quantale homomorphism $\delta$ from
$\hbox{\rm Cen}(B)$ to $\hbox{\rm Cen}(A)$ defined by
$$
\delta(g)(\bigvee_{i}{}_A\hs{x_i\diamond s_i}{y_i}) =
\bigvee_{i}{}_{A}\hs{x_i\diamond g(s_i)}{y_i}
$$ 
for $g\in \hbox{\rm Cen}(B)$ and $x_i, y_i\in X$, $s_i\in B$.
One can easily check that $\delta$ is the inverse of $\gamma$
i.e.\ $\gamma$ provides an involutive quantale isomorphism from
$\hbox{\rm Cen}(A)$ to $\hbox{\rm Cen}(B)$.

(2). 
Let us assume that $A$ and $B$ are commutative i.e.\ each of these 
involutive quantales can be viewed as a complete ${}^{*}$-ideal of its 
centre. 
It is enough to show that $\gamma(A)=B$ and $\delta(B)=A$.
Assume that $r\in A$ and $x, y\in X$. 
Then
$\gamma(r)(\rip{B}{x}{y})=\rip{B}{x}{r\bullet y}$.
Since
$r=\bigvee_{j}\lip{A}{z_{j}\diamond s_j}{t_j}$ 
for suitable $z_j, t_j\in X$ and $s_j\in B$ and $B$ is commutative we have
$$
\begin{array}{lll}
\gamma(r)
(\rip{B}{x}{y})&
=&
\bigvee_{j}\rip{B}{x}{\lip{A}{z_j\diamond s_j}{t_j}\bullet y}\\
&
=&
\bigvee_{j}\rip{B}{x}{(z_j\diamond s_j)\rip{B}{t_j}{y}}\\
&
=&
\bigvee_{j}\rip{B}{x}{z_j\diamond (s_j\cdot \rip{B}{t_j}{y})}\\
&
=&
\bigvee_{j}\rip{B}{x}{z_j\diamond (\rip{B}{t_j}{y}\cdot s_j)}\\
&
=&
\bigvee_{j}\rip{B}{x}{\lip{A}{z_j}{t_j}\bullet{y}}\cdot s_j\\
&
=&
\left(\bigvee_{j} s_j\cdot
\gamma(\lip{A}{z_j}{t_j})\right)\bullet\rip{B}{x}{{y}}.
\end{array}
$$
Since $B$ is an ideal in $\hbox{\rm Cen}(B)$ we have that
$$q = \bigvee_{j} s_j\cdot \gamma(\lip{A}{z_j}{t_j})\in B.$$ 
Hence $\gamma(r)=q$ i.e.\ $\gamma(A)\subseteq B$ and therefore also 
$A\subseteq \delta(B)$. 
Similarly $\delta(B)\subseteq A$ and $B\subseteq \gamma(A)$.
\end{proof}

\begin{lemma}\label{cenisom}
Let $A$ be an m-regular involutive quantale. 
Then $\hbox{\rm Cen}(A)$ and $\hbox{ANat}(A)$ are isomorphic as involutive
quantales.
\end{lemma}

\begin{proof} 
Let $\psi\in \hbox{\rm Cen}(A)$ and let 
$M\in \hbox{\rm m-reg}-Hilb_A$. 
Then, since $A\cong \mathcal{K}_{A}(A)$ as involutive quantales $M$ is a
right Hilbert $\mathcal{K}_{A}(A)$-module.
In particular, we have a right action 
${\diamond}_{{\mathcal{K}}_{A}(A)}: M\times \mathcal{K}_{A}(A)\to M$ 
that gives rise to a non-degenerate involutive quantale homomorphism 
$f:\mathcal{K}_{A}(A)\to \mathcal{A}_{A}(M)$.
Then we can find its unique extension 
$g:\mathcal{A}_{A}(A)\to \mathcal{A}_{A}(M)$
by Corollary 1.5 in \cite{pasekamulti} i.e.\
$M\in MOD_{\mathcal{A}_{A}(A)}$ with the right action 
${\diamond}_{{\mathcal{A}}_{A}(A)} : M\times \mathcal{A}_{A}(A)\to M$.
Moreover the map $\sigma_{M}:M\to M$ given by 
$m\mapsto m{\diamond}_{{\mathcal{A}}_{A}(A)} \psi$
is adjointable. 
Namely,
$$
\begin{array}{lll}
b\cdot \hs{n}{\sigma_{M}(m)}\cdot a&
=&
b\cdot \hs{n}{m\diamond \psi}\cdot a\\
&
=&
b\cdot \hs{n}{(m\diamond \psi)\bullet a}\\
&
=&
b\cdot \hs{n}{m\diamond \psi(a)}\\
&
=&
b\cdot \hs{n}{m}\cdot \psi(a)\\
&
=&
b\cdot \psi(\hs{n}{m}\cdot a)\\
&
=&
\psi(b\cdot \hs{n}{m})\cdot a\\
&
=&
\psi(b)\cdot \hs{n}{m}\cdot a\\
&
=&
\hs{n\diamond \psi(b)^{*}}{m}\cdot a\\
&
=&
\hs{n\diamond \psi^{*}(b^{*})}{m}\cdot a\\
&
=&
\hs{n\diamond \psi^{*}\circ b^{*}}{m}\cdot a\\
&
=&
b\cdot \hs{n\diamond \psi^{*}}{m}\cdot a
.
\end{array}
$$ 
Therefore 
$\hs{n}{\sigma_{M}(m)}=\hs{n\diamond \psi^{*}}{m}$
i.e.\ $\sigma_{M}$ is adjointable.
Hence the map $\Phi:\hbox{\rm Cen}(A)\to \hbox{ANat}(A)$ given by
$\psi\longmapsto (\sigma_{M}:M\to M$, 
$m\mapsto m{\diamond}_{{\mathcal{A}}_{A}(A)} \psi)$ defines
an adjointable natural transformation since any adjointable map preserves
the right action. 
This is easily checked since 
$$
\begin{array}{lll}
(h\circ \sigma_M)(m)\diamond a&
=&
h(\sigma_M(m)\diamond a)\\
&
=&
h(m\diamond \psi \diamond a)\\
&
=&
h(m\diamond \psi(a))\\
&
=&
h(m)\diamond \psi(a)\\
&
=&
h(m)\diamond \psi \diamond a\\
&
=&
\sigma_M(h(m))\diamond a\\
&
=&
(\sigma_M \circ h)(m)\diamond a
\end{array}
$$
i.e.\ $h\circ \sigma_M=\sigma_M \circ h$.
Conversely, if 
$\sigma:\hbox{\rm Id}_{\hbox{\small\rm mreg}-Hilb_A}\to 
\hbox{\rm Id}_{\hbox{\small\rm mreg}-Hilb_A}$ 
is an adjointable natural transformation then $\sigma_A\in \hbox{Cen}(A)$.
\end{proof}

The following theorem provides a fully categorical proof of the first part
of the Theorem \ref{aracoriimp}.

\begin{theorem}\label{coriimp}
Let $A$, $B$ be Morita equivalent m-regular involutive quantales. 
Then $\hbox{ANat}(A)$ and $\hbox{ANat}(B)$ are isomorphic as involutive
quantales.
\end{theorem}

\begin{proof}
Let 
$F: \hbox{\rm mreg}-{Hilb}_A \rightarrow \hbox{\rm mreg}-{Hilb}_B$
be the unitary equivalence functor. 
Then there is a functor
$G: \hbox{\rm mreg}-{Hilb}_B \rightarrow \hbox{\rm mreg}-{Hilb}_A$ 
such that we have natural unitary isomorphisms
$\eta:GF\to \hbox{\rm Id}_{\hbox{\rm\small mreg}-Hilb_A}$ 
and
$\zeta:FG\to \hbox{\rm Id}_{\hbox{\rm\small mreg}-Hilb_B}$.
Similarly as in \cite{gronbak} let us define maps
$T:\hbox{ANat}(A)\to \hbox{ANat}(B)$ and
$S:\hbox{ANat}(B)\to \hbox{ANat}(A)$ 
by the prescription 
$$
\sigma\mapsto (\zeta_{N}\circ F(\sigma_{G(N)})\circ \zeta^{-1}_{N})_{N}\
\mbox{and}\
\rho\mapsto (\eta_{M}\circ G(\rho_{F(M)})\circ \eta^{-1}_{M})_{M}.
$$
Note that since $\eta_{M}$ is a unitary isomorphism we have that
$$
\sigma_{GF(M)}=GF(\sigma_M)\ 
\mbox{and}\
\sigma_{GFG(N)}=GF(\sigma_{G(N)}).$$
So we have the following pair of commuting diagrams, the second diagram is
an application of the functor $F$ to the first diagram.
$$
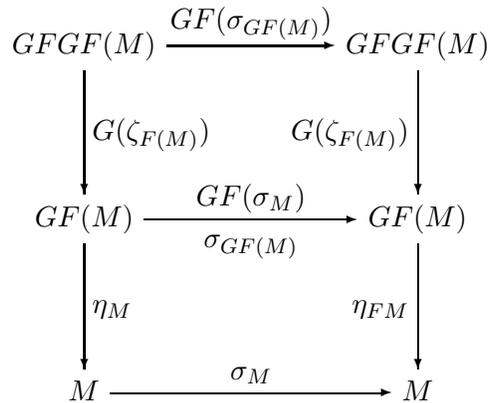
\begin{diagram}
GFGF(M)&
\rTo^{GF(\sigma_{GF(M)})}&
GFGF(M)\\
\dTo_{G(\zeta_{F(M)})}&
&
\dTo^{G(\zeta_{F(M)})}\\
GF(M)&
\rTo^{GF(\sigma_{M})}_{\sigma_{GF(M)}}&
GF(M)\\
\dTo_{\eta_{M}}&
&
\dTo^{\eta_{FM}}\\
M&
\rTo^{\sigma_{M}}&
M\\
\end{diagram}
$$
$$
\begin{diagram}
FGFGF(M)&
\rTo^{FGF(\sigma_{FGF(M)})}&
FGFGF(M)\\
\dTo_{FG(\zeta_{F(M)})}&
&
\dTo^{FG(\zeta_{F(M)})}\\
FGF(M)&
\rTo^{FGF(\sigma_{M})}_{F(\sigma_{GF(M)})}&
FGF(M)\\
\dTo_{F(\eta_{M})}&
&
\dTo^{F(\eta_{FM})}\\
F(M)&
\rTo^{F(\sigma_{M})}&
F(M)\\
\end{diagram}
$$
Hence, from the second diagram, we have
$$F(\sigma_{M})\circ F(\eta_{M})\circ FG(\zeta_{F(M)}) = 
F(\eta_{M})\circ FG(\zeta_{F(M)})\circ FGF(\sigma_{FGF(M)}).$$
In particular,
$$
F(\sigma_{M}) =
F(\eta_{M}) \circ 
FG(\zeta_{F(M)}) \circ 
FGF(\sigma_{FGF(M)}) \circ 
FG(\zeta^{-1}_{F(M)}) \circ 
F(\eta^{-1}_{M})
$$
i.e.\
$$
F(\sigma_{M}) =
F(\eta_{M}) \circ 
FG(T(\sigma)_{F(M)}) \circ 
F(\eta^{-1}_{M})
.
$$
Hence $\sigma_M=ST(\sigma)_M$. 
Similarly, $\rho_N=TS(\rho)_N$ i.e.\
$ST=\hbox{\rm id}_{{\hbox{\rm \small ANat}(A)}}$ and
$TS=\hbox{\rm id}_{{\hbox{\rm \small ANat}(B)}}$. 
Note that we have the following commuting diagram
$$
\begin{diagram}
GFG(N)&
\rTo^{GF(\sigma_{G(N)})}_{\sigma_{GFG(N)}}&
GFG(N)\\
\dTo^{G(\zeta_N)}&
&
\dTo^{G(\zeta_N)}\\
G(N)&
\rTo^{\sigma_{G(N)}}_{G(\zeta_N)\circ 
GF(\sigma_{G(N)})\circ 
G(\zeta_N^{-1})}&
G(N)\\
\end{diagram}
$$
This gives us that 
$$
\sigma_{G(N)} = 
G(\zeta_N)\circ GF(\sigma_{G(N)})\circ G(\zeta_N^{-1})
$$ 
i.e.\ the diagram
$$
\begin{diagram}
G(N_1)&
\rTo^{\sigma_{G(N_1)}}_{G(\zeta_{N_1})\circ 
GF(\sigma_{G(N_{1})})\circ 
G(\zeta^{-1}_{N_{1}})}&
G(N_1)\\
\dTo^{G(f)}&
&\dTo^{G(f)}\\
G(N_2)&
\rTo^{\sigma_{G(N_2)}}_{G(\zeta_{N_2})\circ 
GF(\sigma_{G(N_{2})})\circ 
G(\zeta^{-1}_{N_{2}})}&
G(N_{2})\\
\end{diagram}
$$
commutes. 
Hence also the following diagram commutes
$$
\begin{diagram}
N_1&
\rTo_{\zeta_{N_1}\circ 
F(\sigma_{G(N_{1})})\circ 
\zeta_{N_{1}}^{-1}}&
N_1\\
\dTo^{f}&
&
\dTo^{f}\\
N_2&
\rTo_{\zeta_{N_2}\circ 
F(\sigma_{G(N_{2})})\circ 
\zeta_{N_{2}}^{-1}}&
N_2\\
\end{diagram}
$$
i.e.\ $T(\sigma)$ is a natural transformation and it is evidently 
adjointable. 
Similarly for $S(\rho)$. 
Evidently, $T$ and $S$ are involutive quantale isomorphisms.
\end{proof}

\begin{corollary}\label{unitcom} 
A unital involutive quantale is Morita equivalent to a commutative
m-regular involutive quantale $C$ if and only if it is equivalent to its
own centre.
\end{corollary}

\providecommand{\bysame}{\leavevmode\hbox to3em{\hrulefill}\thinspace}
\providecommand{\MR}{\relax\ifhmode\unskip\space\fi MR }
\providecommand{\MRhref}[2]{%
  \href{http://www.ams.org/mathscinet-getitem?mr=#1}{#2}
}
\providecommand{\href}[2]{#2}

\end{document}